\documentclass[11pt]{article}
\usepackage{amsmath,amsfonts,amssymb, amsthm,enumerate,graphicx}
\usepackage{tikz,authblk}
\usepackage{subfig}
\usepackage{url}
\usepackage{hyperref}
\usepackage{stackengine,scalerel}
\usetikzlibrary {decorations.pathmorphing, decorations.pathreplacing, decorations.shapes}

\usetikzlibrary{arrows.meta}

\newtheorem{theorem} {{\textsf{Theorem}}}
\newtheorem{proposition}[theorem]{{\textsf{Proposition}}}

\newtheorem{definition}[theorem]{{\textsf{Definition}}}

\newtheorem{lemma}[theorem]{{\textsf{Lemma}}}
\newtheorem{question}[theorem]{{\textsf{Question}}}

\newcommand{\xdownarrow}[1]{
	{\left\downarrow\vbox to #1{}\right.\kern-\nulldelimiterspace}
}
\begin{document}

\title{A classification of semi-equivelar gems on the double torus} 
\author{Anshu Agarwal, Biplab Basak$^1$, and Debolina Ghosh}
	
\date{}
	
\maketitle
	
\vspace{-10mm}
\begin{center}
		
\noindent {\small Department of Mathematics, Indian Institute of Technology Delhi, New Delhi 110016, India.$^2$}

\footnotetext[1]{Corresponding author}
		
\footnotetext[2]{{\em E-mail addresses:} \url{maz228084@maths.iitd.ac.in} (A. Agarwal), \url{biplab@iitd.ac.in} (B. Basak), and \url{maz258290@maths.iitd.ac.in} (D. Ghosh).}
		
\medskip
		
\date{August 27, 2026}
\end{center}
	
\hrule
	
\begin{abstract}
A \emph{semi-equivelar gem} of a PL $d$-manifold is a regular colored graph that represents the manifold and admits a regular embedding on a surface, such that the cyclic sequence of face degrees around each vertex is identical. In \cite{ab25, bb24}, semi-equivelar gems of PL $d$-manifolds embedded on surfaces with Euler characteristic $\chi \geq -1$ were classified for $d\geq 2$. 

In this paper, we extend this classification to semi-equivelar gems embedded on the double torus. We show that any such gem must belong to one of the following thirty two types: 
$(4^6)$, $(4^5)$, $(6^4)$, $(4^3,6)$, $(4^3,8)$, $(4^3,12)$, $(4^2,6^2)$, $(4,6,4,6)$, $(4^2,8^2)$, $(4,8,4,8)$, $(8^3)$, $(10^3)$, $(6^2,8)$, $(6^2,10)$, $(6^2,12)$, $(6^2,18)$, $(10^2,4)$, $(12^2,4)$, $(16^2,4)$, $(8^2,6)$, $(12^2,6)$, $(4,6,14)$, $(4,6,16)$, $(4,6,18)$, $(4,6,20)$, $(4,6,24)$, $(4,6,36)$, $(4,8,10)$, $(4,8,12)$, $(4,8,16)$, $(4,8,24)$, and $(4,10,20)$. 
Furthermore, we provide explicit constructions of semi-equivelar gems realizing each of these types.
\end{abstract}


\noindent {\small {\em MSC 2020\,:} Primary 57Q15; Secondary 05C15, 05C10, 52B70, 52C20. 
		
\noindent {\em Keywords:} Semi-equivelar maps; Graph encoded manifold; Semi-equivelar gems; Regular embedding.}
	
\medskip

\section{Introduction}

The theory of graph-encoded manifolds (GEMs) provides a powerful combinatorial framework for the representation and study of piecewise-linear (PL) manifolds. The basic idea is to encode a topological object, namely a PL manifold, by a finite graph endowed with a suitable edge-coloring. This correspondence makes it possible to translate topological questions into finite combinatorial problems and has proved particularly useful in the study of triangulations, minimality, and topological invariants of manifolds.

Let $\Delta_d=\{0,1,\ldots,d\}$. A graph $\Gamma$ is called $(d+1)$-regular if every vertex of $\Gamma$ has degree $d+1$. A $(d+1)$-regular graph $\Gamma$ is said to be a $(d+1)$-regular colored graph if its edges are colored by the elements of $\Delta_d$ in such a way that exactly one edge of each color is incident to every vertex. Given a $(d+1)$-colored graph $(\Gamma,\gamma)$, one can associate a $d$-dimensional simplicial cell complex $\mathcal{K}(\Gamma)$ as follows. For each vertex $u\in V(\Gamma)$, take a $d$-simplex $\sigma(u)$ whose vertices are labeled by the elements of $\Delta_d$. If two vertices $u,v\in V(\Gamma)$ are joined by an edge of color $j$, then the $(d-1)$-faces of $\sigma(u)$ and $\sigma(v)$ opposite to their $j$-labeled vertices are identified so that vertices with the same labels coincide. If $\mathcal{K}(\Gamma)$ is a PL $d$-manifold, then $\Gamma$ is said to encode that manifold. Conversely, a suitable triangulation of a PL manifold gives rise to a $(d+1)$-regular colored graph. Thus, one has the fundamental correspondence between certain types of colored graphs and PL manifolds.

The significance of this correspondence lies in the fact that a topological object can be represented by a relatively simple finite combinatorial object. There is, however, an important distinction between a gem and a crystallization. A gem is a colored graph encoding a manifold, whereas a crystallization is a particularly economical gem. More precisely, a crystallization is a contracted colored graph, in the sense that deleting all the edges of any one color leaves a connected graph. Thus, schematically, the collection of crystallizations of a manifold is contained in the collection of all gems of that manifold.

Crystallizations are particularly useful because they often represent manifolds using a relatively small number of simplices. Consequently, the transition from arbitrary triangulations to colored graphs is especially effective in problems concerning minimal triangulations. Instead of searching among arbitrary triangulations, one can investigate suitably colored graphs and, in particular, crystallizations. Moreover, a gem or crystallization of a PL manifold $M$ provides a combinatorial framework for computing the fundamental group and homology groups of $M$, as well as for studying various other topological properties of $M$.

Suppose that the smallest crystallization of a closed connected PL $d$-manifold $M$ has $2p$ vertices. The \emph{gem-complexity} of $M$ is defined by $k(M)=p-1.$ Since each vertex of a colored graph corresponds to a $d$-simplex in the associated pseudo-complex, minimizing the order of a crystallization is actually minimizing the number of simplices in a triangulation. This provides one of the principal motivations for the use of crystallization theory in the study of minimal triangulations. More generally, an important feature of the gem approach is that topological information about a manifold can be translated into combinatorial information about a finite colored graph.

Another fundamental invariant arising from crystallization theory is the \emph{regular genus}. A $(d+1)$-regular colored graph $\Gamma$ is said to admit a regular embedding into a closed surface $S$ if it can be embedded in $S$ in such a way that every face of the embedding is bounded by a bicolored cycle involving two consecutive colors $\varepsilon_i$ and $\varepsilon_{i+1}$, where the indices are taken modulo $d+1$. The regular genus $\mathcal{G}(M)$ of a closed connected orientable (respectively, non-orientable) PL manifold $M$ is the minimum genus (respectively, half the genus) of an orientable (respectively, non-orientable)  surface into which a gem representing $M$ admits a regular embedding.

The invariants gem-complexity and regular genus provide quantitative information about the topological complexity of manifolds. In particular, they can be related to classical topological invariants which also arise naturally in geometric and variational problems. For example, for a closed connected PL 4-manifold $M$, Basak and Casali~\cite{bc17} established the inequalities
\[
k(M)\geq
3\chi(M)+10\,\operatorname{rk}\pi_1(M)-6 \,\text{ and } \,
\mathcal{G}(M)\geq
2\chi(M)+5\,\operatorname{rk}\pi_1(M)-4.
\]
These inequalities illustrate how purely combinatorial information encoded in a finite colored graph can control fundamental topological quantities of the underlying manifold. Such connections provide a natural motivation for investigating graph-based models in settings where topological and geometric structures interact with variational or nonlinear analytic problems. In particular, finite colored graphs may provide useful combinatorial models for discretizations of differential operators and variational problems on manifolds, including graph-based approximations of operators such as the Laplacian and the $p$-Laplacian. The present work, however, is primarily concerned with the combinatorial and topological structure of the corresponding colored graphs.

The notion of \emph{equivelar} structures was first introduced by McMullen, Schulz and  Wills \cite{Mc82}, where the authors established the existence of equivelar polyhedral manifolds of certain prescribed types. This initiated a systematic study of highly regular decompositions of manifolds. In a subsequent work~\cite{Mc83}, they constructed an infinite family of equivelar $2$-manifolds embedded in $\mathbb{E}^3$, all of which have arbitrarily large genus. Since then, equivelar decompositions of manifolds have been explored extensively from combinatorial, geometric, and topological perspectives.

Nearly two decades later, researchers began investigating the classification of equivelar maps on surfaces of lower genus. In~\cite{dn01}, Datta and Nilakantan achieved a complete classification of all simplicial equivelar polyhedra with up to eleven vertices. Another important contribution was made in~\cite{bk08}, where the authors derived an explicit formula for determining the number of distinct equivelar triangulations and quadrangulations of the torus.

To accommodate broader classes of cell decompositions, the concept of \emph{semi-equivelar} maps was introduced in~\cite{mtu14}. Semi-equivelar maps allow face-cycles at each vertex to follow the same cyclic pattern without requiring all faces to be of the same type, thus generalizing equivelar structures. This framework led to several classification results: in~\cite{dm18}, it was shown that all semi-equivelar maps on the torus and the Klein bottle correspond to Archimedean tilings; in~\cite{dm22}, the authors provided a complete classification of all semi-equivelar maps on the $2$-sphere; and in~\cite{BU121}, a classification was obtained for surfaces of Euler characteristic $-1$. In \cite{c92}, the author introduces the notion of locally regular colored graphs as a higher-dimensional generalization of regular maps on surfaces and describes such graphs on spherical $3$-manifolds.

More recently, the concept has been extended to the setting of colored graphs. The notion of \emph{semi-equivelar gems}  was introduced in~\cite{bb24}. These structures are particularly intriguing: they consist of properly edge-colored graphs embedded on surfaces such that the induced facial cycles form a semi-equivelar map, while the underlying colored graph simultaneously represents a PL $d$-manifold. Thus, semi-equivelar gems provide a rich interplay between combinatorial topology, manifold theory, and the geometry of surface embeddings.

 It is well known that every closed PL $d$-manifold can be represented by a gem, and that a given manifold may admit multiple non-isomorphic gems. A classification of surfaces using gems was established in \cite{b17}. Moreover, any gem admits a regular embedding on a surface \cite{ga81i}. When such an embedding has the property that the cyclic sequence of face degrees around each vertex is identical, the gem is said to be \emph{semi-equivelar}. 

In \cite{bb24}, Basak and Binjola classified all semi-equivelar gems embedded on surfaces with non-negative Euler characteristic. Subsequently, in \cite{ab25}, Agarwal and Basak extended this classification to surfaces with Euler characteristic $-1$. They showed that there exist $15$ possible types of semi-equivelar colored graphs that can embed on such surfaces: one $5$-colored graph, two $4$-colored graphs, and twelve $3$-colored graphs. Furthermore, they proved that a semi-equivelar gem embedded on a surface with Euler characteristic $-1$ cannot be $5$-colored or $4$-colored. For each of the remaining twelve embedding types, they explicitly constructed corresponding $3$-colored semi-equivelar gems.

This article focuses on classifying semi-equivelar gems embedded on the double torus. We first determine all possible types of semi-equivelar colored graphs that can be embedded on this surface (see Lemma~\ref{lemma:possible}). In total, we identify $32$ types: one $6$-colored graph, one $5$-colored graph, eight $4$-colored graphs, and twenty-two $3$-colored graphs. Moreover, we explicitly construct semi-equivelar gems realizing each of these types on the double torus (see Theorems~\ref{thm:construction1}, \ref{thm:construction2},  \ref{thm:construction3}, and \ref{thm:construction4}). 

Hence, the present work provides a complete classification of semi-equivelar gems admitting regular embeddings on the double torus and extends the results of~\cite{ab25, bb24} to the next higher-genus surface. Beyond the classification itself, the results further illustrate how the topology of the embedding surface imposes strong combinatorial restrictions on colored graphs encoding PL manifolds. In this way, semi-equivelar gems form a natural bridge between the combinatorial structure of colored graphs and the global topology of manifolds and their regular embeddings. Such finite combinatorial models are also relevant to the broader topological framework underlying nonlinear analysis, where topological invariants and discrete structures often play an important role in the study of variational problems and their discretizations. In particular, the present classification provides a systematic class of finite graph models of PL manifolds that may serve as a combinatorial framework for future investigations of discrete nonlinear operators and variational problems on manifolds.

\section{Preliminaries}
The theory of edge-colored graphs offers a systematic method for encoding piecewise-linear (PL) manifolds. In particular, every closed, connected PL \(d\)-manifold can be represented by a loopless \((d+1)\)-regular colored graph. This correspondence provides a powerful framework for translating geometric and topological properties of manifolds into combinatorial features of graphs, thereby enabling their study through these graphical models.

\subsection{Graph encoded manifolds (gem)} \label{crystal}
	
Let $\Gamma = (V(\Gamma), E(\Gamma))$ be a multigraph without loops, whose edges are labeled (or colored) by the set $\Delta_d = \{0,1,\ldots,d\}$. An edge-coloring is said to be a \emph{proper edge-coloring} if any two adjacent edges receive distinct colors. The elements of $\Delta_d$ are referred to as the \emph{colors} of $\Gamma$. More precisely, a proper edge-coloring corresponds to a surjective map $\gamma : E(\Gamma) \to \Delta_d$ such that $\gamma(e_1) \neq \gamma(e_2)$ whenever $e_1$ and $e_2$ are adjacent edges. A graph equipped with such a coloring is denoted by $(\Gamma,\gamma)$. A graph $\Gamma$ is called \emph{$(d+1)$-regular} if every vertex has degree $(d+1)$. For standard graph-theoretic terminology, we refer the reader to \cite{bm08}. Throughout this discussion, all spaces are assumed to lie in the PL category.

A \emph{$(d+1)$-regular colored graph} is a pair $(\Gamma,\gamma)$, where $\Gamma$ is a $(d+1)$-regular graph and $\gamma$ is a proper edge-coloring on $\Gamma$. When the coloring is clear from context, we shall simply write $\Gamma$ in place of $(\Gamma,\gamma)$. Given any $(d+1)$-regular colored graph $(\Gamma,\gamma)$, one can associate a $d$-dimensional simplicial cell complex ${\mathcal K}(\Gamma)$, constructed as follows:

\begin{itemize}
\item{} for each vertex $u\in V(\Gamma)$, take a $d$-simplex $\sigma(u)$ with vertices labeled by $\Delta_d$;
		
\item{} corresponding to each edge of color $j$ between $u,v\in V(\Gamma)$, identify the ($d-1$)-faces of $\sigma(u)$ and $\sigma(v)$ opposite to $j$-labeled vertices such that the same labeled vertices coincide.
\end{itemize}
If the geometric carrier $|\mathcal{K}(\Gamma)|$ is (PL) homeomorphic to a PL $d$-manifold $M$, then $\mathcal{K}(\Gamma)$ is called a \emph{colored triangulation} of $M$, and $(\Gamma,\gamma)$ is referred to as a \emph{gem} (graph-encoded manifold) of $M$, or said to \emph{represent} $M$. In particular, every $3$-regular colored graph encodes a closed, connected surface. It is also well known that every closed, connected PL $d$-manifold admits a gem. Moreover, a gem $\Gamma$ representing a manifold $M$ is bipartite if and only if $M$ is orientable.


Now consider a $4$-regular colored graph $(\Gamma,\gamma)$ with color set $\{0,1,2,3\}$. For any distinct $i, j, k \in \{0,1,2,3\}$, let $g_{ij}$ denote the number of connected components of the subgraph $\Gamma_{\{i,j\}}$, and let $g_{ijk}$ denote the number of connected components of the subgraph $\Gamma_{\{i,j,k\}}$.

\begin{proposition}\label{proposition:3-manifold}
Let $(\Gamma,\gamma)$ be a $4$-regular colored graph with $p$ vertices. Then $(\Gamma,\gamma)$ is a gem of a closed connected $3$-manifold if and only if 
$g_{ij}+g_{ik}+g_{jk}=2 g_{ijk}+p/2$ for every $\{i,j,k\}\subset\{0,1,2,3\}$.
\end{proposition}

\subsection{Regular embedding} 

Let $\Gamma$ be a $(d+1)$-regular colored graph. We say that $\Gamma$ embeds regularly on a surface $S$ if it can be embedded on $S$ in such a way that each face of the embedding is bounded by a bi-colored cycle, where the cycle uses two consecutive colors, $\varepsilon_i$ and $\varepsilon_{i+1}$ for some $i$, and indices are taken modulo $d+1$. Here, $\varepsilon = (\varepsilon_0, \dots, \varepsilon_d)$ represents a cyclic permutation of $\Delta_d$.
Regular embeddings are of significant interest in combinatorial topology, and numerous notable results on this topic can be found in \cite{b19, bb21, bc17, cp90, fg82, g81, ga81i}. Below, we present several key results from \cite{ga81i} that are relevant to regular embeddings and will be useful for our article.

\begin{proposition}[\cite{ga81i}]\label{proposition:embedding}
If $\Gamma$ is a bipartite (respectively, non-bipartite) $(d+1)$-regular colored graph representing a closed connected orientable (respectively, non-orientable) PL $d$-manifold $M$, then for each cyclic permutation $\varepsilon = (\varepsilon_0, \varepsilon_1$, $\dots, \varepsilon_d)$ of $\Delta_d$, there exists a regular embedding of $\Gamma$ into an orientable (respectively, non-orientable) surface $S$. 
\end{proposition}


\begin{proposition}[\cite{ga81i}]\label{prop: surface embedding}
A $3$-regular colored graph $\Gamma$ represents a closed connected surface $S$ if and only if it embeds regularly on the surface $S$.   
\end{proposition}



Let $\Gamma$ be a $(d+1)$-regular colored graph embedded regularly on a surface $S$. Since each face in the embedding is bounded by a bi-colored cycle, every face forms a polygon with an even number of sides (including the possibility of a 2-gon). In this article, we restrict our focus to polygons of length at least 4. Without loss of generality, we assume the color sequence $\varepsilon$ to be $(0, 1, \dots, d)$. We define the face-cycles $P_0, P_1, \dots, P_d$ at a vertex $x$ in the embedding of $\Gamma$ on $S$ as the consecutive faces incident to $x$, where each polygon $P_i$ is bounded by a bi-colored cycle of colors $i$ and $i+1$, for $0 \leq i \leq d$, with the condition that $d+1 = 0$.
	   



\begin{definition}\label{def:semiequivelar}
{\rm
Let $\Gamma$ be a $(d+1)$-regular colored graph regularly embedded on a surface $S$. If the face-cycles $P_0,P_1,\ldots,P_d$ incident to each vertex are of the same type, then $\Gamma$ is called a \emph{semi-equivelar graph regularly embedded on $S$}. If $\Gamma$ represents a $d$-manifold $M$, it is called a \emph{semi-equivelar gem of $M$}.}

{\rm
If $P_i$ is a $p_i$-gon for $0\leq i\leq d$, and $\Gamma$ has $p$ vertices, then $\Gamma$ is of type $[(p_0,p_1,\ldots,p_d); p]$. More generally, if each face-cycle contains $n_i$ adjacent $p'_i$-gons, we write the type as $[((p'_0)^{n_0},(p'_1)^{n_1},\ldots$, $(p'_m)^{n_m}); p]$. Here, $p'_i\neq p'_{i+1}$ for $0\leq i<m$, while $p'_i$ may equal $p'_j$ for $|i-j|\geq2$, but $p'_0\neq p'_m$. Since $p$ is determined by the type and $S$, we simply refer to $\Gamma$ as a \emph{semi-equivelar graph of type $((p'_0)^{n_0},\ldots,(p'_m)^{n_m})$ regularly embedded on $S$}. If $\Gamma$ represents a manifold $M$, it is called a \emph{semi-equivelar gem of type $((p'_0)^{n_0},\ldots,(p'_m)^{n_m})$ representing $M$, regularly embedded on $S$}.
}

\end{definition}

\section{Main results}
In~\cite{ab25, bb24}, the authors classified all semi-equivelar gems of PL $d$-manifolds, for $d\geq 2$, whose underlying embeddings lie on surfaces with Euler characteristic $\chi \geq -1$. These results provide a comprehensive understanding of semi-equivelar structures on the sphere, projective plane, torus, Klein bottle, and the surface of Euler characteristic $-1$. 

In the present article, we extend this line of investigation to the next natural case in terms of topological complexity: the double torus, i.e., the closed orientable surface of genus $2$. The following lemma enumerates all possible face-cycle types that can occur for semi-equivelar gems embedded on the double torus. This provides the essential combinatorial groundwork needed for the classification presented in this work.

\begin{lemma}\label{lemma:possible}
If $\Gamma$ is a semi-equivelar graph embedded regularly on a surface $S$ with $\chi(S)=-2$, then $\Gamma$ is one of the following thirty two types: $(4^6)$, $(4^5)$, $(6^4)$, $(4^3,6)$, $(4^3,8)$, $(4^3,12)$, $(4^2,6^2)$, $(4,6,4,6)$, $(4^2,8^2)$, $(4,8,4,8)$, $(8^3)$, $(10^3)$, $(6^2,8)$, $(6^2,10)$, $(6^2,12)$, $(6^2,18)$, $(10^2,4)$, $(12^2,4)$, $(16^2,4)$, $(8^2,6)$, $(12^2,6)$, $(4,6,14)$, $(4,6,16)$, $(4,6,18)$, $(4,6,20)$, $(4,6,24)$, $(4,6,36)$, $(4,8,10)$, $(4,8,12)$, $(4,8,16)$, $(4,8,24)$, and $(4,10,20)$. 
\end{lemma}

\begin{proof}
Let $\Gamma$ be a $[({p_0},{p_1},\dots,{p_d});p]$-type semi-equivelar graph embedded regularly on $S$, where $p_i \geq 4$. Since $\Gamma$ is a regular $(d+1)$-colored graph, $p_i$ is even. Let $q_0,q_1, \dots, q_l$ be the lengths of the polygons of different sizes, where $q_j=p_i$, for some $0\leq i \leq d$. Let $k_i$ be the number of $q_i$-polygons. Thus, $\sum_{i=0}^l{k_i}=d+1$. Let $V, E$, and $F$ denote the number of vertices, edges, and faces in the regular embedding of $\Gamma$ on $S$, respectively. Then $V=p$, $E=p (d+1)/2$, and $F=p(\frac{k_0}{q_0}+\frac{k_1}{q_1}+ \dots +\frac{k_l}{q_l})$. Thus, we have
\begin{align}
\Big(1-\frac{(d+1)}{2}+\frac{k_0}{q_0}+\frac{k_1}{q_1}+ \dots +\frac{k_l}{q_l}\Big)=\frac{\chi(S)}{p}. \label{1}
\end{align}
Since $q_i \geq 4$, we have $k_i/q_i \leq k_i/4$, which further implies
\begin{equation} \label{2}
d+1 \leq 4-\frac{4 \chi(S)}{p}=4+\frac{8}{p}. 
\end{equation}

\noindent Since $p\geq 4, \ (d+1)\leq 6.$ 
		
\noindent \textbf {Case 1.} Let $(d+1)=6$. Equation \eqref{2} gives $p\leq4$, hence we get $p=4$. Since $q_i$ divides $p=4$ and $q_i \geq 4$ so only possibility in this case is $k_0=6, q_0=4$, yielding type $[(4^6);4].$

\noindent \textbf{Case 2.} Let \(d+1 = 5\). The possible tuples \((k_0, k_1, \dots, k_l)\) are 
\((5)\), \((4,1)\), \((3,2)\), \((3,1,1)\), \((2,2,1)\), \((2,1,1,1)\), and \((1,1,1,1,1)\), and Equation~\eqref{1} reduces to
\begin{align}\label{3}
\frac{k_0}{q_0}+\frac{k_1}{q_1}+ \dots + \frac{k_l}{q_l} = \frac{3}{2} - \frac{2}{p}.
\end{align}

\textbf{Case \((5)\):} If \(q_0 \ge 6\) then since $q_0$ divides $p$ so we must have $p\geq 6$, then $\frac{3}{2}=\frac{5}{q_0}+\frac{2}{p}\leq \frac{5}{6}+\frac{2}{6}=\frac{7}{6}$, a contradiction. For \(q_0 = 4\), Equation~\eqref{1} gives \(p=8\), yielding type \(\left[(4^{5});\,8\right]\).

\textbf{Case \((4,1)\):} Equation~\eqref{3} becomes \(\frac{4}{q_0} + \frac{1}{q_1} + \frac{2}{p} = \frac{3}{2}\). First, assume \(q_0 > q_1\). Since \(q_1 \ge 4\), we obtain
$
\frac{3}{2}
= \frac{4}{q_0} + \frac{1}{q_1} + \frac{2}{p}
\le \frac{4}{6} + \frac{1}{4} + \frac{2}{6}
= \frac{5}{4}.
$
This is a contradiction; hence \(q_0 > q_1\) is not possible.

Next, assume \(q_0 < q_1\). If \(q_0 \ge 6\), then
$
\frac{3}{2}
= \frac{4}{q_0} + \frac{1}{q_1} + \frac{2}{p}
\le \frac{4}{6} + \frac{1}{8} + \frac{2}{8}
= \frac{25}{24},
$
which again yields a contradiction. Therefore, \(q_0 = 4\).
Now if \(q_1 \ge 8\). Then
$
\frac{3}{2}
= \frac{4}{q_0} + \frac{1}{q_1} + \frac{2}{p}
\le \frac{4}{4} + \frac{1}{8} + \frac{2}{8}
= \frac{11}{8},
$
which is impossible. Thus, the only remaining candidate is \(q_1 = 6\), giving
$
q_0 = 4$, $p = 6.
$
However, \(4 \nmid 6\), so this configuration is not admissible. Hence, for the tuple \((4,1)\), no semi-equivelar graph can be regularly embedded in a surface of Euler characteristic \(-2\).

\textbf{Case \((3,2)\):} Equation~\eqref{3} gives \(\frac{3}{q_0} + \frac{2}{q_1} + \frac{2}{p} = \frac{3}{2}\). First, assume \(q_0 < q_1\). Then
$
\frac{3}{2}
\le \frac{3}{4} + \frac{2}{6} + \frac{2}{6}
= \frac{17}{12}
$
which is a contradiction. Thus, \(q_0 < q_1\) is not possible.
Next, assume \(q_1 < q_0\). Then
$
\frac{3}{2}
\le \frac{3}{6} + \frac{2}{4} + \frac{2}{6}
= \frac{4}{3}
$
which again leads to a contradiction. Hence, \(q_1 < q_0\) is also impossible.
Therefore, for the tuple \((3,2)\), there is no semi-equivelar graph that can be regularly embedded in the surface \(\#_{2}\mathbb{T}^{2}\).

\textbf{Case \((3,1,1)\):} \(\frac{3}{q_0} + \frac{1}{q_1} + \frac{1}{q_2} + \frac{2}{p} = \frac{3}{2}\). Any admissible integral choices of \(p\) and \(q_i\) fail to satisfy the equality; hence, no embedding exists.

\textbf{Case \((2,2,1)\):} \(\frac{2}{q_0} + \frac{2}{q_1} + \frac{1}{q_2} + \frac{2}{p} = \frac{3}{2}\). Any admissible integral choices of \(p\) and \(q_i\) fail to satisfy the equality; hence, no embedding exists.

\textbf{Case \((2,1,1,1)\):} \(\frac{2}{q_0} + \frac{1}{q_1} + \frac{1}{q_2} + \frac{1}{q_3} + \frac{2}{p} = \frac{3}{2}\). Any admissible integral choices of \(p\) and \(q_i\) fail to satisfy the equality; hence, no embedding exists.

\textbf{Case \((1,1,1,1,1)\):} \(\frac{1}{q_0} + \dots + \frac{1}{q_4} + \frac{2}{p} = \frac{3}{2}\). Any admissible integral choices of \(p\) and \(q_i\) fail to satisfy the equality; hence, no embedding exists.

Hence, the only \(5\)-colored semi-equivelar graph regularly embedded on \(\#_{2}\mathbb{T}^{2}\) is of type \(\left[(4^{5});\,8\right]\).

\noindent \textbf{Case 3.} Let \(d+1=4\), with possible tuples \((k_0,\dots,k_l)\) being \((4),(3,1),(2,2),(2,1,1),(1,1,1,1)\). Equation~\eqref{1} becomes
\begin{equation}\label{4}
\frac{k_0}{q_0} + \dots + \frac{k_l}{q_l} = 1 - \frac{2}{p}.
\end{equation}

\textbf{Case \((4)\):} \(\frac{4}{q_0} + \frac{2}{p} = 1\). Only \(q_0=6\) gives \(p=6\), yielding type \(\left[(6^4);6\right]\).

\textbf{Case \((3,1)\):} \(\frac{3}{q_0} + \frac{1}{q_1} + \frac{2}{p} = 1\). Considering \(q_0>q_1\) or \(q_0<q_1\) and admissible ranges, we get
\[
(q_0,q_1,p) = (4,6,24),\ (4,8,16),\ (4,12,12),
\]
giving types \([(4^3,6);24], [(4^3,8);16], [(4^3,12);12]\).

\textbf{Case \((2,2)\):} \(\frac{2}{q_0} + \frac{2}{q_1} + \frac{2}{p} = 1\). Only feasible pairs are \((q_0,q_1)=(4,6),(4,8)\), giving types \([(4^2,6^2);12]\), \([(4^2,8^2);8]\), $[(4,6,4,6);12]$, and  $[(4,8,4,8);8]$.

\textbf{Case \((2,1,1)\):} \(\frac{2}{q_0} + \frac{1}{q_1} + \frac{1}{q_2} + \frac{2}{p} = 1\). Any admissible integral choices of \(p\) and \(q_i\) fail to satisfy the equality; hence, no embedding exists.

\textbf{Case \((1,1,1,1)\):} \(\frac{1}{q_0} + \frac{1}{q_1} + \frac{1}{q_2} + \frac{1}{q_3} + \frac{2}{p} = 1\). Any admissible integral choices of \(p\) and \(q_i\) fail to satisfy the equality; hence, no embedding exists.

\noindent\textbf{Case 4.} Let \(d+1=3\). The possible tuples are \((k_0,\dots,k_l) = (3), (2,1), (1,1,1)\), and Equation~\eqref{1} becomes
\begin{equation}\label{9}
\frac{k_0}{q_0} + \dots + \frac{k_l}{q_l} = \frac{1}{2} - \frac{2}{p}.
\end{equation}

\textbf{Case \((3)\):} \(\frac{3}{q_0} + \frac{2}{p} = \frac{1}{2}\). Only \(q_0=8,10\) yield integral \(p\), giving types \([(8^3);16]\) and \([(10^3);10]\).

\textbf{Case \((2,1)\):} \(\frac{2}{q_0} + \frac{1}{q_1} + \frac{2}{p} = \frac{1}{2}\).  
Analyzing \(q_0>q_1\) and \(q_0<q_1\) gives admissible types:
\[
[(6^2,8);48],\ [(6^2,10);30],\ [(6^2,12);24],\ [(6^2,18);18],\ 
[(10^2,4);40],\ [(12^2,4);24],\]
\[ [(16^2,4);16],\ [(8^2,6);24],\ [(12^2,6);12].
\]

\textbf{Case \((1,1,1)\):} \(\frac{1}{q_0} + \frac{1}{q_1} + \frac{1}{q_2} + \frac{2}{p} = \frac{1}{2}\).  
Only \(q_0=4\) gives solutions. The admissible \((q_1,q_2)\) pairs are
\[
(6,14),(6,16),(6,18),(6,20),(6,24),(6,36),\ 
(8,10),(8,12),(8,16),(8,24),(10,20),
\]
yielding embedding types
\[
[(4,6,14);168],\ [(4,6,16);96],\ [(4,6,18);72],\ [(4,6,20);60],\ [(4,6,24);48],\ [(4,6,36);36], 
\]
\[
[(4,8,10);80],\ [(4,8,12);48],\ [(4,8,16);32],\ [(4,8,24);24],\ [(4,10,20);20].
\]

Thus, all 3-colored semi-equivelar graphs regularly embedded in \(\#_2\mathbb{T}^2\) are of types:
\[
(8^3),\ (10^3),\ (6^2,8),\ (6^2,10),\ (6^2,12),\ (6^2,18),\ (10^2,4),\ (12^2,4),\ (16^2,4),\ (8^2,6),\ (12^2,6),
\]
\[
(4,6,14),\ (4,6,16),\ (4,6,18),\ (4,6,20),\ (4,6,24),\ (4,6,36),\ (4,8,10),\ (4,8,12),\ (4,8,16),\]
\[(4,8,24),\ (4,10,20).
\]
This completes the proof.
\end{proof}

For each possible type of \(\Gamma\) listed in Lemma~\ref{lemma:possible}, we now demonstrate the existence of a regular embedding on the double torus.

\begin{theorem}\label{thm:construction1}
There exists a semi-equivelar gem of type $(4^6)$ that is regularly embedded on the double torus.
\end{theorem}

\begin{proof}
Let $(\Gamma,\gamma)$ be a semi-equivelar gem of type $(4^6)$. Then the total number of vertices is $4$. Since each component of any 4-colored subgraph of $(\Gamma,\gamma)$ represents $\mathbb{S}^3$, by Proposition~\ref{proposition:3-manifold}, for every $\{i,j,k\}\subset\{0, 1, 2, 3, 4, 5\}$ we have
\[
g_{ij} + g_{ik} + g_{jk} = 2 g_{ijk} + 2,
\]
where $g_{ij}$ denotes the number of connected components of the subgraph $\Gamma_{\{i,j\}}$, and $g_{ijk}$ denotes the number of connected components of the subgraph $\Gamma_{\{i, j, k\}}$.
\begin{figure}[h]
   \begin{center}
\tikzset{every picture/.style={line width=0.75pt}} 


\caption{Embedding on $\#_2 \mathbb{T}^2$ of a gem representing $ \mathbb{S}^5$ of type $(4^6)$.} 
    \label{fig:0}
   \end{center}
\end{figure}

 One admissible choice of the parameters $g_{ij}$ and $g_{ijk}$ for $i,j,k\in\{0, 1, 2, 3, 4, 5\}$ is the following: $g_{02}=g_{04}=g_{13}=g_{15}=g_{24}=g_{35}=2$, all remaining $g_{ij}=1$, $g_{024}=g_{135}=2$, and the remaining $g_{ijk}=1$.
In Figure~\ref{fig:0}, we construct a semi-equivelar graph of type $(4^6)$ realizing these parameters. The subgraph induced on any four color is connected. Since it is a gem,  the known characterization of $3$-dimensional gems with $4$ vertices implies that the represented manifold is $\mathbb{S}^3$. Now, add one further color to the $4$-colored subgraph. The resulting $5$-colored subgraph represents a $4$-dimensional colored complex whose vertex links are precisely the $3$-dimensional complexes represented by the corresponding $4$-colored subgraphs. Hence, the link of every vertex is $\mathbb{S}^3$, and therefore the resulting complex represents a closed $4$-dimensional combinatorial manifold. Since the corresponding graph has $4$ vertices, the known characterization of $4$-dimensional gems with $4$ vertices implies that the represented manifold is $\mathbb{S}^4$. Adding the remaining color, we obtain the full $6$-colored graph $(\Gamma,\gamma)$. The link of every vertex in the corresponding $5$-dimensional complex is $\mathbb{S}^4$. Hence, the complex is a closed $5$-dimensional combinatorial manifold. Since $(\Gamma,\gamma)$ has exactly $4$ vertices, the known characterization of $5$-dimensional gems with $4$ vertices implies that the represented manifold is $\mathbb{S}^5$.
 Therefore, Figure~\ref{fig:0} gives a semi-equivelar gem of type $(4^6)$ regularly embedded on a surface $S$ representing $\mathbb{S}^5$. The surface $S$ has three vertices, labeled $a$, $b$, and $c$. The vertex $a$ lies in the interior of the region bounded by the $(2,3)$-colored $4$-cycle, the vertex $b$ lies in the interior of the region bounded by the $(3,4)$-colored $4$-cycle and the vertex $c$ lies in the interior of the region bounded by the $(4,5)$-colored $4$-cycle. The \(1\)-cells, depicted as dotted segments on the boundary, are identified by matching the corresponding \(x_i\)'s. Since $S$ is orientable, and has three vertices, six edges and one face, it follows that $S$ is a double torus. Thus, Figure~\ref{fig:0} depicts a semi-equivelar gem of type $(4^6)$ representing $\mathbb{S}^4$, regularly embedded on the surface $\#_2 \mathbb{T}^2$.
\end{proof}

\begin{theorem}\label{thm:construction2}
There exists a semi-equivelar gem of type $(4^5)$ that is regularly embedded on the double torus.
\end{theorem}
\begin{proof}
Let $(\Gamma,\gamma)$ be a semi-equivelar gem of type $(4^5)$. Then the total number of vertices is $8$. Since each component of any 4-colored subgraph of $(\Gamma,\gamma)$ represents $\mathbb{S}^3$, by Proposition~\ref{proposition:3-manifold}, for every $\{i,j,k\}\subset\{0,1,2,3,4\}$ we have
\[
g_{ij} + g_{ik} + g_{jk} = 2 g_{ijk} + 4,
\]
where $g_{ij}$ denotes the number of connected components of the subgraph $\Gamma_{\{i,j\}}$, and $g_{ijk}$ denotes the number of connected components of the subgraph $\Gamma_{\{i,j,k\}}$.

\begin{figure}[h]
\begin{center}
\tikzset{every picture/.style={line width=0.75pt}} 



\caption{Embedding on $\#_2 \mathbb{T}^2$ of a gem representing $ \mathbb{S}^4$ of type $(4^5)$.}
    \label{fig:1}
    
\end{center}
\end{figure}

 One admissible choice of the parameters $g_{ij}$ and $g_{ijk}$ for $i,j,k\in\{0,1,2,3, 4\}$ is the following: $g_{03}=g_{14}=4$, all remaining $g_{ij}=2$, $g_{012}=g_{024}=g_{123}=g_{234}=1$, and the remaining $g_{ijk}=2$.
In Figure~\ref{fig:1}, we construct a semi-equivelar graph of type $(4^5)$ realizing these parameters. Moreover, it is easy to check that the $3$-manifold represented by each component of any 4-colored subgraph of $(\Gamma,\gamma)$ is indeed $\mathbb{S}^3$. Therefore, Figure~\ref{fig:1} gives a semi-equivelar gem of type $(4^5)$ regularly embedded on a surface $S$. The surface $S$ has four vertices, labeled $a$, $b$, $c$, and $d$. Each of the vertices $a$ and $b$ lies in the interior of the region bounded by a $(0,1)$-colored $4$-cycle, while each of the vertices $c$ and $d$ lies in the interior of the region bounded by a $(1,2)$-colored $4$-cycle.
 The \(1\)-cells, depicted as dotted segments on the boundary, are identified by matching the corresponding \(x_i\)'s. Since $S$ is orientable and has four vertices, seven edges and one face, it follows that $S$ is a double torus. Finally, it is straightforward to check that $(\Gamma,\gamma)$ represents $\mathbb{S}^4$. Thus, Figure~\ref{fig:1} depicts a semi-equivelar gem of type $(4^5)$ representing $\mathbb{S}^4$, regularly embedded on the surface $\#_2 \mathbb{T}^2$.
\end{proof}

\begin{theorem}\label{thm:construction3}
For each of the following eight types: $(6^4)$,  $(4^3,6)$, $(4^3,8)$, $(4^3,12)$, $(4^2,6^2)$, $(4^2,8^2)$,  $(4,6,4,6)$, and $(4,8,4,8)$, there exists a semi-equivelar gem that is regularly embedded on the double torus.
\end{theorem}

\begin{proof}
Let $(\Gamma,\gamma)$ be a semi-equivelar gem of one of the following eight types: 
$(6^4)$, $(4^3,6)$, $(4^3,8)$, $(4^3,12)$, $(4^2,6^2)$, $(4^2,8^2)$, $(4,6,4,6)$, and $(4,8,4,8)$. 
Let $p$ denote the number of vertices of $(\Gamma,\gamma)$. Then, by Proposition~\ref{proposition:3-manifold}, for every $\{i,j,k\}\subset\{0,1,2,3\}$ we have
\[
g_{ij} + g_{ik} + g_{jk} = 2 g_{ijk} + \frac{p}{2},
\]
where $g_{ij}$ denotes the number of connected components of the subgraph $\Gamma_{\{i,j\}}$, and $g_{ijk}$ denotes the number of connected components of the subgraph $\Gamma_{\{i,j,k\}}$.

Figures~\ref{fig:2}--\ref{fig:6} illustrate a CW-complex structure of a surface with a single \(2\)-cell, on which the graph \(\Gamma\) is regularly embedded. The boundary of the \(2\)-cell contains at most five vertices from \(\{a,b,c,d,e\}\), and the \(1\)-cells, depicted as dotted segments on the boundary, are identified by matching the corresponding \(x_i\)'s. We now describe each figure in detail.

\smallskip

\noindent {\bf \underline{$(6^4)$-type}}: For the type $(6^4)$, the number of vertices is $6$. As stated in Proposition~\ref{proposition:3-manifold}, one possible choice of the parameters $g_{ij}$ and $g_{ijk}$ for $i,j,k \in \{0,1,2,3\}$ is:
\[
g_{01}=g_{12}=g_{23}=g_{03}=1,\qquad 
g_{13}=g_{02}=3,\qquad 
g_{023}=g_{013}=g_{012}=g_{123}=1,\qquad 
\]

In Figure~\ref{fig:2}, we construct a semi-equivelar gem of type $(6^4)$ with $6$ vertices representing $\mathbb{S}^3$, regularly embedded on a surface $S$. The surface $S$ has three vertices, labeled $a$, $b$, and $c$, where $a$ is lying in the interior of a region bounded by a $(2,3)$-colored $6$-cycle, $b$ is lying in the interior of a region bounded by a $(1,2)$-colored $6$-cycle and $c$ is lying in the interior of a region bounded by a $(0,3)$-colored $6$-cycle. Since $S$ is orientable and has six edges and one face, it follows that $S$ is a double torus. Thus, Figure~\ref{fig:2} depicts a semi-equivelar gem of type $(6^4)$ regularly embedded on the surface $\#_2 \mathbb{T}^2$.

\smallskip

\begin{figure}[h]
    \begin{center}

\tikzset{every picture/.style={line width=0.44pt}} 

\begin{tikzpicture}[x=0.45pt,y=0.45pt,yscale=-1,xscale=1]

\draw   (369.33,271.67) -- (325.5,347.59) -- (237.83,347.59) -- (194,271.67) -- (237.83,195.75) -- (325.5,195.75) -- cycle ;
\draw    (325.5,195.75) -- (369.33,271.67) (333.96,202.41) -- (327.04,206.41)(338.96,211.07) -- (332.04,215.07)(343.96,219.73) -- (337.04,223.73)(348.96,228.39) -- (342.04,232.39)(353.96,237.05) -- (347.04,241.05)(358.96,245.71) -- (352.04,249.71)(363.96,254.37) -- (357.04,258.37)(368.96,263.03) -- (362.04,267.03) ;
\draw    (237.83,347.59) -- (325.5,347.59) (247.83,343.59) -- (247.83,351.59)(257.83,343.59) -- (257.83,351.59)(267.83,343.59) -- (267.83,351.59)(277.83,343.59) -- (277.83,351.59)(287.83,343.59) -- (287.83,351.59)(297.83,343.59) -- (297.83,351.59)(307.83,343.59) -- (307.83,351.59)(317.83,343.59) -- (317.83,351.59) ;
\draw    (194,271.67) -- (237.83,195.75) (195.54,261.01) -- (202.46,265.01)(200.54,252.35) -- (207.46,256.35)(205.54,243.69) -- (212.46,247.69)(210.54,235.03) -- (217.46,239.03)(215.54,226.37) -- (222.46,230.37)(220.54,217.71) -- (227.46,221.71)(225.54,209.04) -- (232.46,213.04)(230.54,200.38) -- (237.46,204.38) ;
\draw  [dash pattern={on 0.84pt off 2.51pt}] (140.17,147.67) -- (423.17,147.67) -- (423.17,395.67) -- (140.17,395.67) -- cycle ;
\draw  [dash pattern={on 4.5pt off 4.5pt}]  (245.33,146) -- (237.83,195.75) ;
\draw  [dash pattern={on 4.5pt off 4.5pt}]  (314.33,147) -- (325.5,195.75) ;
\draw  [dash pattern={on 4.5pt off 4.5pt}]  (423.33,311) -- (369.33,271.67) ;
\draw  [dash pattern={on 4.5pt off 4.5pt}]  (139.33,318) -- (194,271.67) ;
\draw  [dash pattern={on 4.5pt off 4.5pt}]  (183.17,393.92) -- (237.83,347.59) ;
\draw  [dash pattern={on 4.5pt off 4.5pt}]  (325.5,347.59) -- (384.33,395) ;
\draw    (182.33,149) .. controls (184.68,148.8) and (185.96,149.87) .. (186.16,152.22) .. controls (186.36,154.57) and (187.63,155.64) .. (189.98,155.44) .. controls (192.33,155.24) and (193.6,156.31) .. (193.81,158.66) .. controls (194.01,161.01) and (195.28,162.08) .. (197.63,161.88) .. controls (199.98,161.69) and (201.25,162.76) .. (201.45,165.11) .. controls (201.66,167.46) and (202.93,168.53) .. (205.28,168.33) .. controls (207.63,168.13) and (208.9,169.2) .. (209.1,171.55) .. controls (209.31,173.9) and (210.58,174.97) .. (212.93,174.77) .. controls (215.28,174.57) and (216.55,175.64) .. (216.75,177.99) .. controls (216.96,180.34) and (218.23,181.41) .. (220.58,181.21) .. controls (222.93,181.01) and (224.2,182.08) .. (224.4,184.43) .. controls (224.6,186.78) and (225.87,187.85) .. (228.22,187.65) .. controls (230.57,187.45) and (231.84,188.52) .. (232.05,190.87) .. controls (232.25,193.22) and (233.52,194.29) .. (235.87,194.09) -- (237.83,195.75) -- (237.83,195.75) ;
\draw    (377.33,147) .. controls (377.26,149.36) and (376.05,150.5) .. (373.69,150.43) .. controls (371.34,150.36) and (370.12,151.5) .. (370.05,153.85) .. controls (369.98,156.2) and (368.76,157.35) .. (366.41,157.28) .. controls (364.06,157.21) and (362.84,158.35) .. (362.76,160.7) .. controls (362.69,163.05) and (361.47,164.2) .. (359.12,164.13) .. controls (356.77,164.06) and (355.55,165.2) .. (355.48,167.55) .. controls (355.41,169.9) and (354.19,171.05) .. (351.84,170.98) .. controls (349.49,170.91) and (348.27,172.05) .. (348.19,174.4) .. controls (348.12,176.75) and (346.9,177.9) .. (344.55,177.83) .. controls (342.2,177.76) and (340.98,178.9) .. (340.91,181.25) .. controls (340.84,183.6) and (339.62,184.75) .. (337.27,184.68) .. controls (334.92,184.61) and (333.7,185.75) .. (333.62,188.1) .. controls (333.55,190.45) and (332.33,191.6) .. (329.98,191.53) .. controls (327.63,191.46) and (326.41,192.61) .. (326.34,194.96) -- (325.5,195.75) -- (325.5,195.75) ;
\draw    (424.33,219) .. controls (424.28,221.36) and (423.08,222.51) .. (420.72,222.46) .. controls (418.37,222.41) and (417.16,223.57) .. (417.11,225.92) .. controls (417.06,228.27) and (415.85,229.42) .. (413.5,229.37) .. controls (411.15,229.32) and (409.94,230.48) .. (409.89,232.83) .. controls (409.84,235.18) and (408.63,236.34) .. (406.28,236.29) .. controls (403.93,236.24) and (402.72,237.4) .. (402.67,239.75) .. controls (402.62,242.11) and (401.41,243.26) .. (399.05,243.21) .. controls (396.7,243.16) and (395.49,244.31) .. (395.44,246.66) .. controls (395.39,249.01) and (394.18,250.17) .. (391.83,250.12) .. controls (389.48,250.07) and (388.27,251.23) .. (388.22,253.58) .. controls (388.17,255.93) and (386.96,257.09) .. (384.61,257.04) .. controls (382.26,256.99) and (381.05,258.15) .. (381,260.5) .. controls (380.95,262.85) and (379.74,264.01) .. (377.39,263.96) .. controls (375.04,263.91) and (373.83,265.06) .. (373.78,267.41) .. controls (373.73,269.77) and (372.52,270.92) .. (370.16,270.87) -- (369.33,271.67) -- (369.33,271.67) ;
\draw    (138.33,219) .. controls (140.69,218.93) and (141.9,220.08) .. (141.97,222.44) .. controls (142.04,224.79) and (143.25,225.94) .. (145.6,225.87) .. controls (147.95,225.81) and (149.16,226.96) .. (149.23,229.31) .. controls (149.3,231.66) and (150.51,232.81) .. (152.86,232.75) .. controls (155.21,232.68) and (156.42,233.83) .. (156.49,236.18) .. controls (156.56,238.54) and (157.77,239.69) .. (160.13,239.62) .. controls (162.48,239.55) and (163.69,240.7) .. (163.76,243.05) .. controls (163.83,245.4) and (165.04,246.55) .. (167.39,246.49) .. controls (169.74,246.43) and (170.95,247.58) .. (171.02,249.93) .. controls (171.09,252.28) and (172.3,253.43) .. (174.65,253.36) .. controls (177.01,253.29) and (178.22,254.44) .. (178.29,256.8) .. controls (178.36,259.15) and (179.57,260.3) .. (181.92,260.24) .. controls (184.27,260.17) and (185.48,261.32) .. (185.55,263.67) .. controls (185.62,266.02) and (186.83,267.17) .. (189.18,267.11) .. controls (191.53,267.04) and (192.74,268.19) .. (192.81,270.54) -- (194,271.67) -- (194,271.67) ;
\draw    (237.83,347.59) .. controls (239.92,348.68) and (240.42,350.27) .. (239.33,352.36) .. controls (238.24,354.45) and (238.74,356.04) .. (240.83,357.13) .. controls (242.92,358.22) and (243.41,359.81) .. (242.32,361.9) .. controls (241.23,363.99) and (241.73,365.58) .. (243.82,366.67) .. controls (245.91,367.76) and (246.41,369.35) .. (245.32,371.44) .. controls (244.23,373.53) and (244.72,375.12) .. (246.81,376.21) .. controls (248.9,377.3) and (249.4,378.89) .. (248.31,380.98) .. controls (247.22,383.07) and (247.72,384.66) .. (249.81,385.75) .. controls (251.9,386.85) and (252.39,388.44) .. (251.3,390.53) .. controls (250.21,392.62) and (250.71,394.21) .. (252.8,395.3) -- (253.33,397) -- (253.33,397) ;
\draw    (325.5,347.59) .. controls (326.79,349.56) and (326.44,351.19) .. (324.47,352.48) .. controls (322.5,353.77) and (322.15,355.4) .. (323.44,357.37) .. controls (324.73,359.34) and (324.39,360.98) .. (322.42,362.27) .. controls (320.45,363.56) and (320.1,365.19) .. (321.39,367.16) .. controls (322.68,369.13) and (322.33,370.76) .. (320.36,372.05) .. controls (318.39,373.34) and (318.04,374.98) .. (319.33,376.95) .. controls (320.62,378.92) and (320.28,380.55) .. (318.31,381.84) .. controls (316.34,383.13) and (315.99,384.76) .. (317.28,386.73) .. controls (318.57,388.7) and (318.22,390.34) .. (316.25,391.63) -- (315.33,396) -- (315.33,396) ;
\draw  [fill={rgb, 255:red, 0; green, 0; blue, 0 }  ,fill opacity=1 ] (138.77,148.84) .. controls (138.81,147.99) and (139.74,147.32) .. (140.86,147.36) .. controls (141.97,147.4) and (142.85,148.13) .. (142.81,148.99) .. controls (142.77,149.85) and (141.84,150.51) .. (140.72,150.47) .. controls (139.61,150.43) and (138.74,149.7) .. (138.77,148.84) -- cycle ;
\draw  [fill={rgb, 255:red, 0; green, 0; blue, 0 }  ,fill opacity=1 ] (214.77,147.84) .. controls (214.81,146.99) and (215.74,146.32) .. (216.86,146.36) .. controls (217.97,146.4) and (218.85,147.13) .. (218.81,147.99) .. controls (218.77,148.85) and (217.84,149.51) .. (216.72,149.47) .. controls (215.61,149.43) and (214.74,148.7) .. (214.77,147.84) -- cycle ;
\draw  [fill={rgb, 255:red, 0; green, 0; blue, 0 }  ,fill opacity=1 ] (275.77,147.84) .. controls (275.81,146.99) and (276.74,146.32) .. (277.86,146.36) .. controls (278.97,146.4) and (279.85,147.13) .. (279.81,147.99) .. controls (279.77,148.85) and (278.84,149.51) .. (277.72,149.47) .. controls (276.61,149.43) and (275.74,148.7) .. (275.77,147.84) -- cycle ;
\draw  [fill={rgb, 255:red, 0; green, 0; blue, 0 }  ,fill opacity=1 ] (338.77,147.84) .. controls (338.81,146.99) and (339.74,146.32) .. (340.86,146.36) .. controls (341.97,146.4) and (342.85,147.13) .. (342.81,147.99) .. controls (342.77,148.85) and (341.84,149.51) .. (340.72,149.47) .. controls (339.61,149.43) and (338.74,148.7) .. (338.77,147.84) -- cycle ;
\draw  [fill={rgb, 255:red, 0; green, 0; blue, 0 }  ,fill opacity=1 ] (420.77,148.84) .. controls (420.81,147.99) and (421.74,147.32) .. (422.86,147.36) .. controls (423.97,147.4) and (424.85,148.13) .. (424.81,148.99) .. controls (424.77,149.85) and (423.84,150.51) .. (422.72,150.47) .. controls (421.61,150.43) and (420.74,149.7) .. (420.77,148.84) -- cycle ;
\draw  [fill={rgb, 255:red, 0; green, 0; blue, 0 }  ,fill opacity=1 ] (420.77,268.84) .. controls (420.81,267.99) and (421.74,267.32) .. (422.86,267.36) .. controls (423.97,267.4) and (424.85,268.13) .. (424.81,268.99) .. controls (424.77,269.85) and (423.84,270.51) .. (422.72,270.47) .. controls (421.61,270.43) and (420.74,269.7) .. (420.77,268.84) -- cycle ;
\draw  [fill={rgb, 255:red, 0; green, 0; blue, 0 }  ,fill opacity=1 ] (419.77,394.84) .. controls (419.81,393.99) and (420.74,393.32) .. (421.86,393.36) .. controls (422.97,393.4) and (423.85,394.13) .. (423.81,394.99) .. controls (423.77,395.85) and (422.84,396.51) .. (421.72,396.47) .. controls (420.61,396.43) and (419.74,395.7) .. (419.77,394.84) -- cycle ;
\draw  [fill={rgb, 255:red, 0; green, 0; blue, 0 }  ,fill opacity=1 ] (347.77,395.84) .. controls (347.81,394.99) and (348.74,394.32) .. (349.86,394.36) .. controls (350.97,394.4) and (351.85,395.13) .. (351.81,395.99) .. controls (351.77,396.85) and (350.84,397.51) .. (349.72,397.47) .. controls (348.61,397.43) and (347.74,396.7) .. (347.77,395.84) -- cycle ;
\draw  [fill={rgb, 255:red, 0; green, 0; blue, 0 }  ,fill opacity=1 ] (281.77,395.84) .. controls (281.81,394.99) and (282.74,394.32) .. (283.86,394.36) .. controls (284.97,394.4) and (285.85,395.13) .. (285.81,395.99) .. controls (285.77,396.85) and (284.84,397.51) .. (283.72,397.47) .. controls (282.61,397.43) and (281.74,396.7) .. (281.77,395.84) -- cycle ;
\draw  [fill={rgb, 255:red, 0; green, 0; blue, 0 }  ,fill opacity=1 ] (216.77,396.84) .. controls (216.81,395.99) and (217.74,395.32) .. (218.86,395.36) .. controls (219.97,395.4) and (220.85,396.13) .. (220.81,396.99) .. controls (220.77,397.85) and (219.84,398.51) .. (218.72,398.47) .. controls (217.61,398.43) and (216.74,397.7) .. (216.77,396.84) -- cycle ;
\draw  [fill={rgb, 255:red, 0; green, 0; blue, 0 }  ,fill opacity=1 ] (138.77,394.84) .. controls (138.81,393.99) and (139.74,393.32) .. (140.86,393.36) .. controls (141.97,393.4) and (142.85,394.13) .. (142.81,394.99) .. controls (142.77,395.85) and (141.84,396.51) .. (140.72,396.47) .. controls (139.61,396.43) and (138.74,395.7) .. (138.77,394.84) -- cycle ;
\draw  [fill={rgb, 255:red, 0; green, 0; blue, 0 }  ,fill opacity=1 ] (137.77,268.84) .. controls (137.81,267.99) and (138.74,267.32) .. (139.86,267.36) .. controls (140.97,267.4) and (141.85,268.13) .. (141.81,268.99) .. controls (141.77,269.85) and (140.84,270.51) .. (139.72,270.47) .. controls (138.61,270.43) and (137.74,269.7) .. (137.77,268.84) -- cycle ;
\draw    (480.82,213) -- (513.33,213) ;
\draw    (480.35,250) -- (517.05,249.33) (490.28,245.82) -- (490.42,253.82)(500.27,245.64) -- (500.42,253.64)(510.27,245.46) -- (510.42,253.45) ;
\draw    (485.44,289) .. controls (487.11,287.33) and (488.77,287.33) .. (490.44,289) .. controls (492.11,290.67) and (493.77,290.67) .. (495.44,289) .. controls (497.11,287.33) and (498.77,287.33) .. (500.44,289) .. controls (502.11,290.67) and (503.77,290.67) .. (505.44,289) .. controls (507.11,287.33) and (508.77,287.33) .. (510.44,289) .. controls (512.11,290.67) and (513.77,290.67) .. (515.44,289) -- (519.73,289) -- (519.73,289) ;
\draw  [dash pattern={on 4.5pt off 4.5pt}]  (485.53,326) -- (519.82,326) ;

\draw (306,127.4) node [anchor=north west][inner sep=0.75pt]  [font=\footnotesize]  {$x_{1}$};
\draw (434,305.4) node [anchor=north west][inner sep=0.75pt]  [font=\footnotesize]  {$x_{1}$};
\draw (381,402.4) node [anchor=north west][inner sep=0.75pt]  [font=\footnotesize]  {$x_{2}$};
\draw (173,403.4) node [anchor=north west][inner sep=0.75pt]  [font=\footnotesize]  {$x_{2}$};
\draw (120,308.4) node [anchor=north west][inner sep=0.75pt]  [font=\footnotesize]  {$x_{3}$};
\draw (240,125.4) node [anchor=north west][inner sep=0.75pt]  [font=\footnotesize]  {$x_{3}$};
\draw (173,126.4) node [anchor=north west][inner sep=0.75pt]  [font=\footnotesize]  {$x_{4}$};
\draw (373,126.4) node [anchor=north west][inner sep=0.75pt]  [font=\footnotesize]  {$x_{4}$};
\draw (433,212.4) node [anchor=north west][inner sep=0.75pt]  [font=\footnotesize]  {$x_{5}$};
\draw (311,399.4) node [anchor=north west][inner sep=0.75pt]  [font=\footnotesize]  {$x_{5}$};
\draw (248,402.4) node [anchor=north west][inner sep=0.75pt]  [font=\footnotesize]  {$x_{6}$};
\draw (119,209.4) node [anchor=north west][inner sep=0.75pt]  [font=\footnotesize]  {$x_{6}$};
\draw (240,201.4) node [anchor=north west][inner sep=0.75pt]  [font=\tiny]  {$v_{0}$};
\draw (311,203.4) node [anchor=north west][inner sep=0.75pt]  [font=\tiny]  {$v_{1}$};
\draw (342,269.4) node [anchor=north west][inner sep=0.75pt]  [font=\tiny]  {$v_{2}$};
\draw (310,326.4) node [anchor=north west][inner sep=0.75pt]  [font=\tiny]  {$v_{3}$};
\draw (238,330.4) node [anchor=north west][inner sep=0.75pt]  [font=\tiny]  {$v_{4}$};
\draw (206,271.4) node [anchor=north west][inner sep=0.75pt]  [font=\tiny]  {$v_{5}$};
\draw (211,123.73) node [anchor=north west][inner sep=0.75pt]  [font=\footnotesize]  {$a$};
\draw (119,259.73) node [anchor=north west][inner sep=0.75pt]  [font=\footnotesize]  {$a$};
\draw (214,406.73) node [anchor=north west][inner sep=0.75pt]  [font=\footnotesize]  {$a$};
\draw (346,404.73) node [anchor=north west][inner sep=0.75pt]  [font=\footnotesize]  {$a$};
\draw (435,262.73) node [anchor=north west][inner sep=0.75pt]  [font=\footnotesize]  {$a$};
\draw (337,130.73) node [anchor=north west][inner sep=0.75pt]  [font=\footnotesize]  {$a$};
\draw (428,134.73) node [anchor=north west][inner sep=0.75pt]  [font=\footnotesize]  {$b$};
\draw (280,405.73) node [anchor=north west][inner sep=0.75pt]  [font=\footnotesize]  {$b$};
\draw (125,129.73) node [anchor=north west][inner sep=0.75pt]  [font=\footnotesize]  {$b$};
\draw (423,406.73) node [anchor=north west][inner sep=0.75pt]  [font=\footnotesize]  {$c$};
\draw (131,397.73) node [anchor=north west][inner sep=0.75pt]  [font=\footnotesize]  {$c$};
\draw (272,126.73) node [anchor=north west][inner sep=0.75pt]  [font=\footnotesize]  {$c$};
\draw (492.62,218.4) node [anchor=north west][inner sep=0.75pt]    {$0$};
\draw (494.27,256.4) node [anchor=north west][inner sep=0.75pt]    {$1$};
\draw (496.56,295.4) node [anchor=north west][inner sep=0.75pt]    {$2$};
\draw (496.33,331.4) node [anchor=north west][inner sep=0.75pt]    {$3$};

\draw   (194, 271.67) circle [x radius= 5, y radius= 5]   ;
\draw   (237.83, 195.75) circle [x radius= 5, y radius= 5]   ;
\draw   (237.83, 195.75) circle [x radius= 5, y radius= 5]   ;
\draw   (194, 271.67) circle [x radius= 5, y radius= 5]   ;
\draw   (237.83, 195.75) circle [x radius= 5, y radius= 5]   ;
\draw   (194, 271.67) circle [x radius= 5, y radius= 5]   ;
\draw   (237.83, 347.59) circle [x radius= 5, y radius= 5]   ;
\draw   (237.83, 347.59) circle [x radius= 5, y radius= 5]   ;
\draw   (325.5, 347.59) circle [x radius= 5, y radius= 5]   ;
\draw   (237.83, 347.59) circle [x radius= 5, y radius= 5]   ;
\draw   (325.5, 347.59) circle [x radius= 5, y radius= 5]   ;
\draw   (369.33, 271.67) circle [x radius= 5, y radius= 5]   ;
\draw   (325.5, 195.75) circle [x radius= 5, y radius= 5]   ;
\draw   (325.5, 195.75) circle [x radius= 5, y radius= 5]   ;
\draw   (369.33, 271.67) circle [x radius= 5, y radius= 5]   ;
\draw   (325.5, 195.75) circle [x radius= 5, y radius= 5]   ;
\draw   (369.33, 271.67) circle [x radius= 5, y radius= 5]   ;
\end{tikzpicture}

\caption{Embedding on $\#_2 \mathbb{T}^2$ of a gem representing $\mathbb{S}^3$ of type $(6^4)$.}
    \label{fig:2}
    \end{center}
\end{figure}
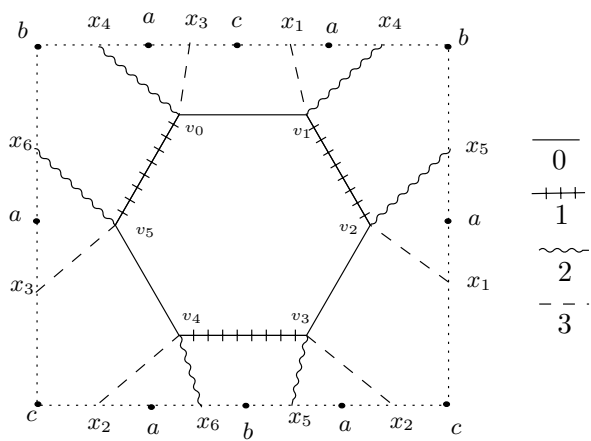

\noindent {\bf \underline{$(4^3,6)$-type}}: For the type \((4^3,6)\), there are 24 vertices. By Proposition~\ref{proposition:3-manifold}, one possible choice of parameters \(g_{ij}\) and \(g_{ijk}\) for \(i,j,k \in \{0,1,2,3\}\) is
\[
g_{01}=g_{12}=g_{23}=6,\quad 
g_{03}=g_{13}=g_{02}=4,\quad 
g_{023}=g_{013}=1,\quad 
g_{012}=g_{123}=2.
\]

Figure~\ref{fig:3} illustrates a semi-equivelar gem of type \((4^3,6)\) with 24 vertices representing \(L(5,2)\), regularly embedded on a surface \(S\). The surface has six vertices \(a,b,c,d,e,f\), each lying inside a region bounded by a \((1,2)\)-colored 4-cycle. Since \(S\) is orientable with nine edges and one face, it is a double torus. Hence, the figure depicts a semi-equivelar gem of type \((4^3,6)\) regularly embedded on \(\#_2 \mathbb{T}^2\).

\smallskip

\begin{figure}[h]
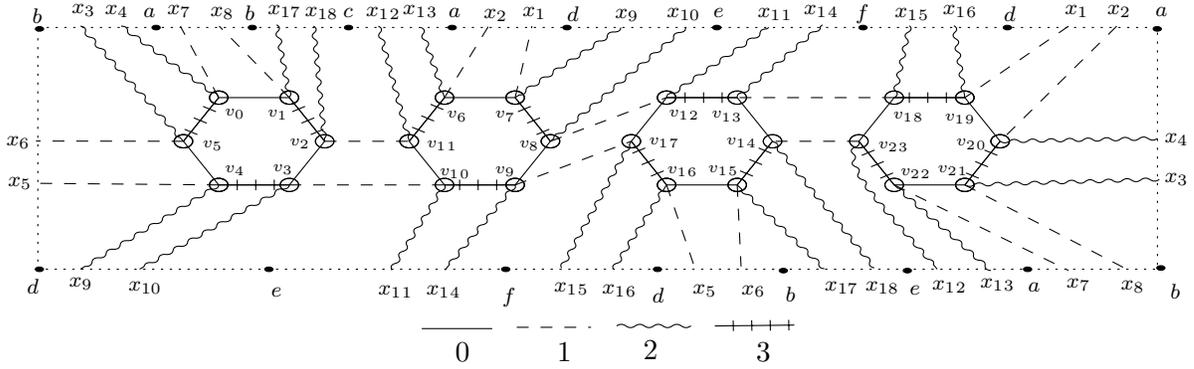

   \begin{center}
\tikzset{every picture/.style={line width=0.4pt}} 


\caption{Embedding on $\#_2 \mathbb{T}^2$ of a gem representing $L(5,2)$ of type $(4^3,6)$.}
    \label{fig:3}
   \end{center}
\end{figure}

\noindent {\bf \underline{$(4^3,8)$-type}}: For the type $(4^3,8)$, the number of vertices is $16$. As stated in Proposition~\ref{proposition:3-manifold}, one possible choice of the parameters $g_{ij}$ and $g_{ijk}$ for $i,j,k \in \{0,1,2,3\}$ is:
\[
g_{01}=g_{12}=g_{23}=g_{13}=g_{02}=4,\qquad 
g_{03}=2,\qquad 
g_{023}=g_{013}=1,\qquad 
g_{012}=g_{123}=2.
\]

In Figure~\ref{fig:4}$(a)$, we construct a semi-equivelar gem of type $(4^3,8)$ with $16$ vertices representing $\mathbb{RP}^3$, regularly embedded in a surface $S$. The surface $S$ has four vertices $a,b,c,d$, each lying in the interior of a region bounded by a $(1,2)$-colored $4$-cycle, and has seven edges and one face; hence $S$ is orientable of genus two. Therefore, Figure~\ref{fig:4}$(a)$ illustrates a semi-equivelar gem of type $(4^3,8)$ regularly embedded on $\#_2 \mathbb{T}^2$

\smallskip

\begin{figure}[h]
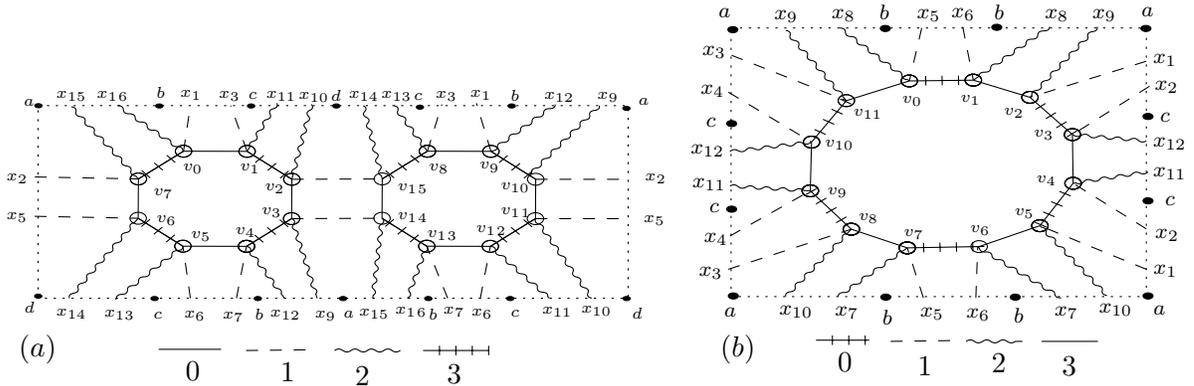

   \begin{center}
\tikzset{every picture/.style={line width=0.4pt}} 


\caption{Embedding on $\#_2 \mathbb{T}^2$: $(a)$ a gem representing $\mathbb{RP}^3$ of type $(4^3,8)$; $(b)$ a gem representing $\mathbb{S}^3$ of type $(4^3,12)$.}
    \label{fig:4}
       \end{center}
\end{figure}

\noindent {\bf \underline{$(4^3,12)$-type}}: For the type \((4^3,12)\), there are 12 vertices. By Proposition~\ref{proposition:3-manifold}, one possible choice of parameters \(g_{ij}\) and \(g_{ijk}\) for \(i,j,k \in \{0,1,2,3\}\) is
\[
g_{01}=g_{12}=g_{23}=3,\quad 
g_{03}=1,\quad
g_{13}=g_{02}=4,\quad
g_{023}=g_{013}=1,\quad 
g_{012}=g_{123}=2.
\]

Figure~\ref{fig:4}$(b)$ shows a semi-equivelar gem of type \((4^3,12)\) with 12 vertices representing \(\mathbb{RP}^3\), regularly embedded on a surface \(S\). The surface has three vertices \(a,b,c\), each inside a region bounded by a \((1,2)\)-colored 4-cycle. Since \(S\) is orientable with six edges and one face, it is a double torus. Hence, the figure depicts a semi-equivelar gem of type \((4^3,12)\) regularly embedded on \(\#_2 \mathbb{T}^2\).

\smallskip

\noindent {\bf \underline{$(4^2,6^2)$-type}}: For the type $(4^2,6^2)$, the number of vertices is $12$. As stated in Proposition~\ref{proposition:3-manifold}, one possible choice of the parameters $g_{ij}$ and $g_{ijk}$ for $i,j,k \in \{0,1,2,3\}$ is:
\[
g_{01}=g_{12}=g_{13}=3,\qquad 
g_{03}=g_{23}=2,\qquad
g_{02}=4,\qquad
g_{023}=g_{013}=g_{123}=1,\qquad 
g_{012}=2.
\]

In Figure~\ref{fig:5}$(a)$, we construct a semi-equivelar gem of type $(4^2,6^2)$ with $12$ vertices representing $\mathbb{RP}^3$, regularly embedded on a surface $S$. The surface $S$ has five vertices, labeled $a$, $b$, $c$, $d$ and $e$, where each of the vertices $a$, $b$, and $c$ lies in the interior of the region bounded by a $(0,1)$-colored $4$-cycle, while each of the vertices $d$ and $e$ lies in the interior of the region bounded by a  $(0,3)$-colored $6$-cycle.
 Since $S$ is orientable and has eight edges and one face, it follows that $S$ is a double torus. Thus, Figure~\ref{fig:5}$(a)$ depicts a semi-equivelar gem of type $(4^3,12)$ regularly embedded on the surface $\#_2 \mathbb{T}^2$.

\smallskip

\begin{figure}[h]
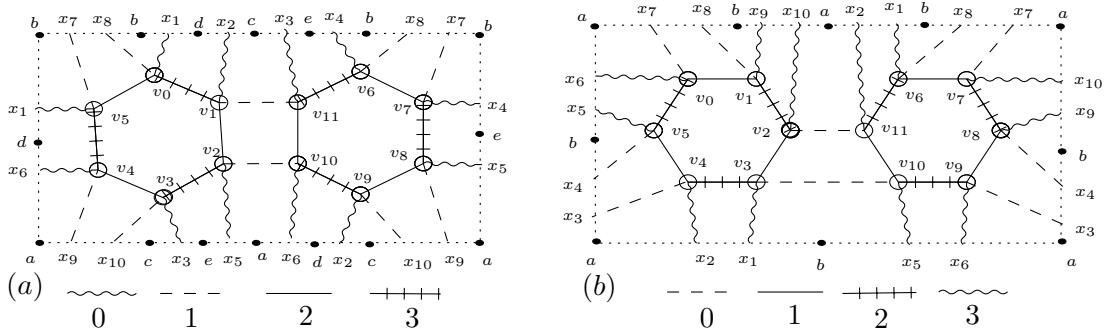

   \begin{center}
\tikzset{every picture/.style={line width=0.4pt}} 



\caption{Embedding on $\#_2 \mathbb{T}^2$: $(a)$ a gem representing $ \mathbb{RP}^3$ of type $(4^2,6^2)$; $(b)$ a gem representing $L(3,1)$ of type $(4,6,4,6)$.}
    \label{fig:5}
    
   \end{center}
\end{figure}

\noindent {\bf \underline{$(4,6,4,6)$-type}}: For the type \((4,6,4,6)\), there are 12 vertices. By Proposition~\ref{proposition:3-manifold}, one possible choice of parameters \(g_{ij}\) and \(g_{ijk}\) for \(i,j,k \in \{0,1,2,3\}\) is
\[
g_{01}=g_{23}=3,\quad 
g_{03}=g_{12}=2,\quad
g_{02}=g_{13}=3,\quad
g_{012}=g_{013}=g_{023}=g_{123}=1.
\]

Figure~\ref{fig:5}$(b)$ illustrates a semi-equivelar gem of type \((4,6,4,6)\) with 12 vertices representing \(L(3,1)\), regularly embedded on a surface \(S\). The surface has two vertices, \(a\) and \(b\); each of the vertices is lying inside a region bounded by a \((0,3)\)-colored 6-cycle. Since \(S\) is orientable with five edges and one face, it is a double torus. Thus, the figure depicts a semi-equivelar gem of type \((4,6,4,6)\) regularly embedded on \(\#_2 \mathbb{T}^2\).

\smallskip

\begin{figure}[h]
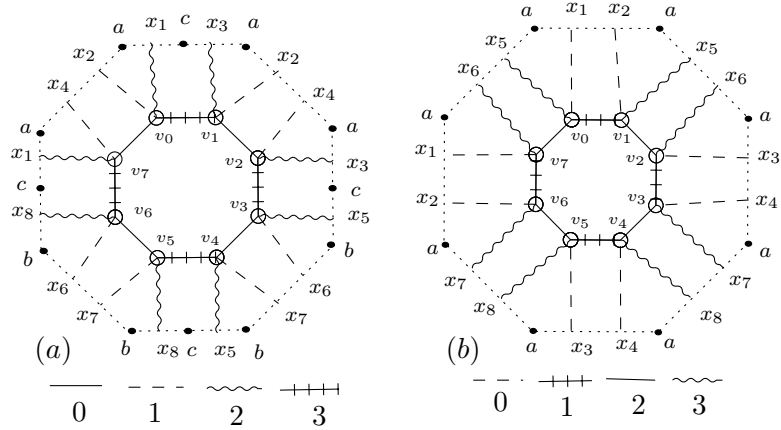

\begin{center}
\tikzset{every picture/.style={line width=0.4pt}} 

\caption{Embedding on $\#_2 \mathbb{T}^2$: $(a)$ a gem representing $ \mathbb{S}^3$ of type $(4^2,8^2)$; $(b)$ a gem representing $ \mathbb{S}^3$ of type $(4,8,4,8)$.}
    \label{fig:6}
\end{center}
\end{figure}

\smallskip

\noindent {\bf \underline{$(4^2,8^2)$-type}}: For the type \((4^2,8^2)\), there are 8 vertices. By Proposition~\ref{proposition:3-manifold}, one possible choice of parameters \(g_{ij}\) and \(g_{ijk}\) for \(i,j,k \in \{0,1,2,3\}\) is
\[
g_{01}=g_{12}=2,\quad 
g_{23}=g_{03}=1,\quad
g_{02}=4, \quad 
g_{13}=3, \quad
g_{012}=2, \quad
g_{013}=g_{023}=g_{123}=1.
\]

Figure~\ref{fig:6}$(a)$ shows a semi-equivelar gem of type \((4^2,8^2)\) with 8 vertices representing \(\mathbb{S}^3\), regularly embedded on a surface \(S\). The surface has three vertices \(a,b,c\), where each of the vertices \(a,b\) lies inside the region bounded by a \((1,2)\)-colored 4-cycle and the vertex \(c\) lies inside a region bounded by a \((2,3)\)-colored 6-cycle. Since \(S\) is orientable with six edges and one face, it is a double torus. Hence, the figure depicts a semi-equivelar gem of type \((4^2,8^2)\) regularly embedded on \(\#_2 \mathbb{T}^2\).

\smallskip

\noindent {\bf \underline{$(4,8,4,8)$-type}}: For the type \((4,8,4,8)\), there are 8 vertices. By Proposition~\ref{proposition:3-manifold}, one possible choice of the parameters \(g_{ij}\) and \(g_{ijk}\) for \(i,j,k \in \{0,1,2,3\}\) is
\[
g_{01}=g_{23}=2,\quad 
g_{12}=g_{03}=1,\quad
g_{02}=g_{13}=3, \quad 
g_{012}=g_{013}=g_{023}=g_{123}=1.
\]

Figure~\ref{fig:6}$(b)$ shows a semi-equivelar gem of type \((4,8,4,8)\) with 8 vertices representing \(\mathbb{S}^3\), regularly embedded on a surface \(S\). The surface has a single vertex \(a\) inside a region bounded by a \((0,3)\)-colored 8-cycle. Since \(S\) is orientable with four edges and one face, it is a double torus. Hence, the figure depicts a semi-equivelar gem of type \((4,8,4,8)\) regularly embedded on \(\#_2 \mathbb{T}^2\).
\end{proof}

\begin{theorem}\label{thm:construction4}
For each of the following twenty two types: $(8^3)$, $(10^3)$, $(6^2,8)$, $(6^2,10)$, $(6^2,12)$, $(6^2,18)$, $(10^2,4)$, $(12^2,4)$, $(16^2,4)$, $(8^2,6)$, $(12^2,6)$, $(4,6,14)$, $(4,6,16)$, $(4,6,18)$, $(4,6,20)$, $(4,6,24)$, $(4,6,36)$, $(4,8,10)$, $(4,8,12)$, $(4,8,16)$, $(4,8,24)$, and $(4,10,20)$, there exists a semi-equivelar gem that is regularly embedded on the double torus.
\end{theorem}

\begin{proof}
It follows from Proposition~\ref{prop: surface embedding} that a 3-regular colored graph \(\Gamma\) represents a closed connected surface \(S\) if and only if it is regularly embedded on \(S\). Consequently, if \(\Gamma\) is regularly embedded on \(\#_{2}\mathbb{T}^{2}\), it represents the surface \(\#_{2}\mathbb{T}^{2}\) itself. 

Figures~\ref{fig:7}--\ref{fig:22} illustrate a CW-complex structure of a surface on which \(\Gamma\) embeds regularly, consisting of a single 2-cell. The boundary of the 2-cell contains at most sixteen vertices from \(\{a,b,c,d,e,f,g,h,i,j,k,l,m,n,o,p\}\). The 1-cells, shown as dotted lines on the boundary, are identified by matching the corresponding \(x_i\)'s. We now provide a detailed description of each figure.

\begin{figure}[h]
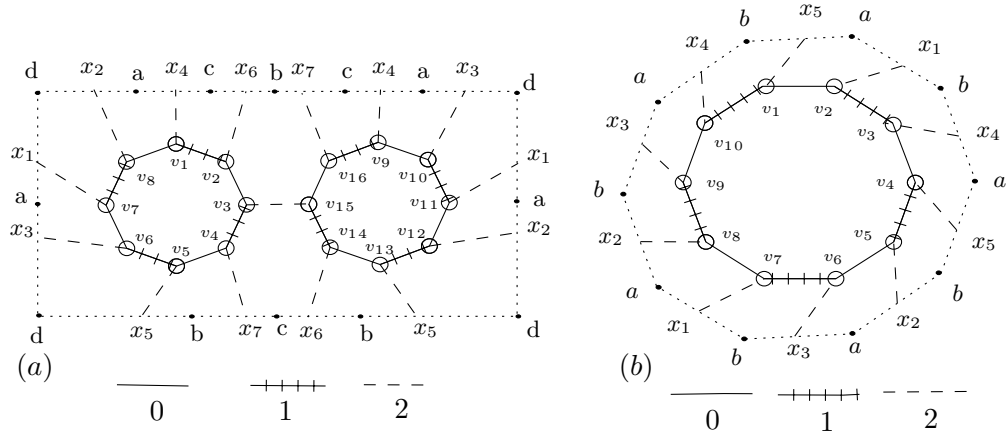

    \begin{center}
\tikzset{every picture/.style={line width=0.4pt}} 



\caption{Embedding on $\#_2 \mathbb{T}^2$ of a gem representing $\#_{2} \mathbb{T}^2$: $(a)$ of type $(8^3)$; $(b)$ of type  $(10^3)$.}
    \label{fig:7}
    \end{center}
\end{figure}

\smallskip

\noindent {\bf \underline{$(8^{3})$-type}}: The type $(8^3)$ consists of 16 vertices.
 In Figure \ref{fig:7}$(a)$, the surface comprises four vertices $a$, $b$, $c$, and $d$, together with seven edges and a single face; since the surface is orientable, it is a double torus. 
The $(0,1)$-colored cycle bounds the two inner $8$-gons 
$v_{1}v_{2}v_{3}v_{4}v_{5}v_{6}v_{7}v_{8}$ and 
$v_{9}v_{10}v_{11}v_{12}v_{13}v_{14}v_{15}v_{16}$. 
The $(0,2)$-colored $8$-cycle encloses the regions 
$v_{1}v_{8}v_{12}v_{11}v_{7}v_{6}v_{10}v_{9}$ and 
$v_{5}v_{4}v_{16}v_{15}v_{3}v_{2}v_{14}v_{13}$, containing the vertices $a$ and $b$, respectively, in their interiors. 
Similarly, the $(1,2)$-colored $8$-cycle bounds the regions 
$v_{1}v_{2}v_{14}v_{15}v_{3}v_{4}v_{16}v_{9}$ and 
$v_{8}v_{7}v_{11}v_{10}v_{6}v_{5}v_{13}v_{12}$, containing the vertices $c$ and $d$, respectively in their interiors. 
Hence, this figure represents a semi-equivelar gem of type $[(8^{3});16]$ regularly embedded on the surface $\#_{2}\mathbb{T}^{2}$.

\smallskip

\noindent {\bf \underline{$(10^{3})$-type}}: For the type $(10^3)$, there are 10 vertices. 
 In Figure \ref{fig:7}$(b)$, the surface is described by the two vertices 
$a$ and $b$, together with five edges and a single face. 
As the surface is orientable, it is a double torus 
$\#_{2}\mathbb{T}^{2}$. The vertex $a$ lies in the interior of the region encircled by the $(0,2)$-colored $10$-cycle, whereas the vertex $b$ occupies the region bounded by the $(1,2)$-colored $10$-cycle. Additionally, the $(0,1)$-colored cycle constitutes the boundary of the central $10$-gon. Consequently, this configuration gives rise to a semi-equivelar gem of type $\left[(10^{3});\,10\right]$ that is regularly embedded in the surface $\#_{2}\mathbb{T}^{2}$.

\smallskip

\begin{figure}[h]
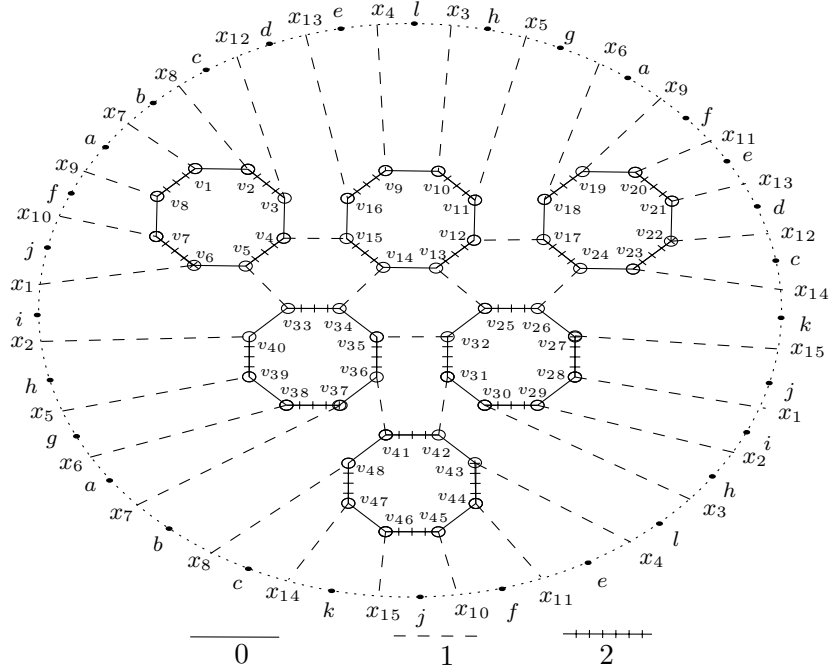

   \begin{center}
       \tikzset{every picture/.style={line width=0.4pt}} 


\caption{Embedding on $\#_2 \mathbb{T}^2$ of a gem representing $\#_{2} \mathbb{T}^2$ of type $(6^2,8)$.}
    \label{fig:8}

   \end{center}
\end{figure}

\noindent {\bf \underline{$(6^{2},\,8)$-type}}: For the type $(6^2,8)$, the corresponding semi-equivelar graph has 48 vertices.
 In Figure \ref{fig:8}, the surface comprises twelve vertices 
$a, b, c, d, e, f, g, h, i, j, k,$ and $l$, together with fifteen edges and a single face. Since the surface is orientable, it is a double torus $\#_{2}\mathbb{T}^{2}$. Each of the vertices $a, c, e, h,$ and $j$ lies in the interior of the region enclosed by a $(1,2)$-colored $6$-cycle, while each of the vertices $b, d, f, g, i, k,$ and $l$ lies in the interior of the region enclosed by a $(0,1)$-colored $6$-cycle. Furthermore, the $(0,2)$-colored cycles bound the innermost six $8$-gons. Accordingly, this configuration yields a semi-equivelar gem of type $\left[(6^{2},\,8);\,48\right]$ regularly embedded in the surface $\#_{2}\mathbb{T}^{2}$.

\smallskip

\noindent {\bf \underline{$(6^2,10)$-type}}: The type $(6^2,10)$ consists of 30 vertices. In Figure \ref{fig:9}$(a)$, the underlying surface contains nine vertices $a$, $b$, $c$, $d$, $e$, $f$, $g$, $h$ and $i$ together with twelve edges and a single face; since the surface is orientable, it is a double torus. 
Each of the vertices $a$, $c$, $e$, $h$ are located in the interior of the region bounded by a $(1,2)$-colored $6$-cycle, whereas each of the vertex $b$, $d$, $f$, $g$, $i$ lies in the interior of the region bounded by a $(0,1)$-colored $6$-cycle. 
The $(0,2)$-colored cycles bound the inner three $10$-gons. 
Consequently, this figure represents a semi-equivelar gem of type $[(6^2,10);30]$ regularly embedded on the surface $\#_{2}\mathbb{T}^{2}$.

\begin{figure}[h]
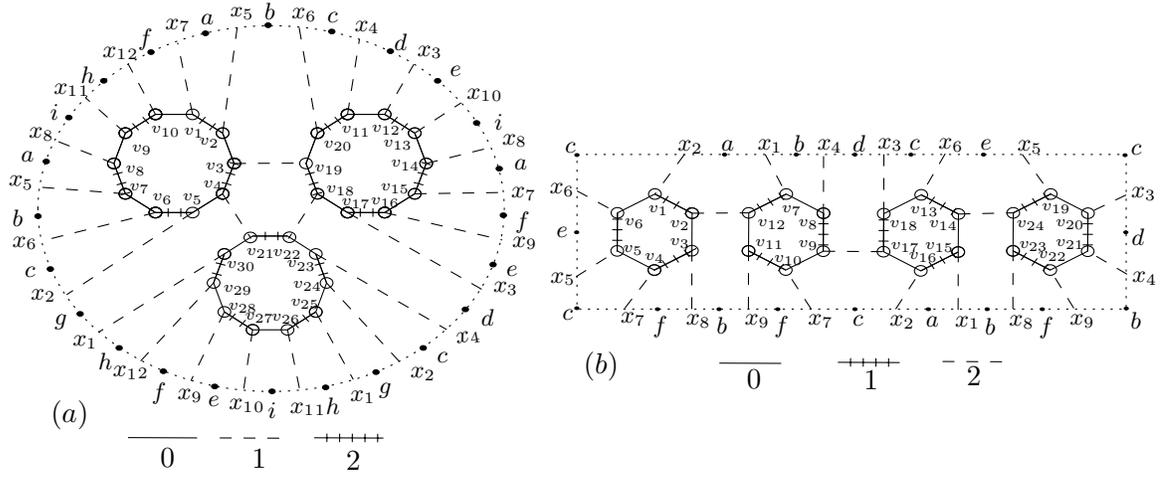

  \begin{center}
\tikzset{every picture/.style={line width=0.4pt}} 



\caption{Embedding on $\#_2 \mathbb{T}^2$ of a gem representing $\#_{2} \mathbb{T}^2$: $(a)$  of type $(6^2,10)$; $(b)$ of type  $(6^2,12)$.}
    \label{fig:9}
 
  \end{center}
\end{figure}

\begin{figure}[h]
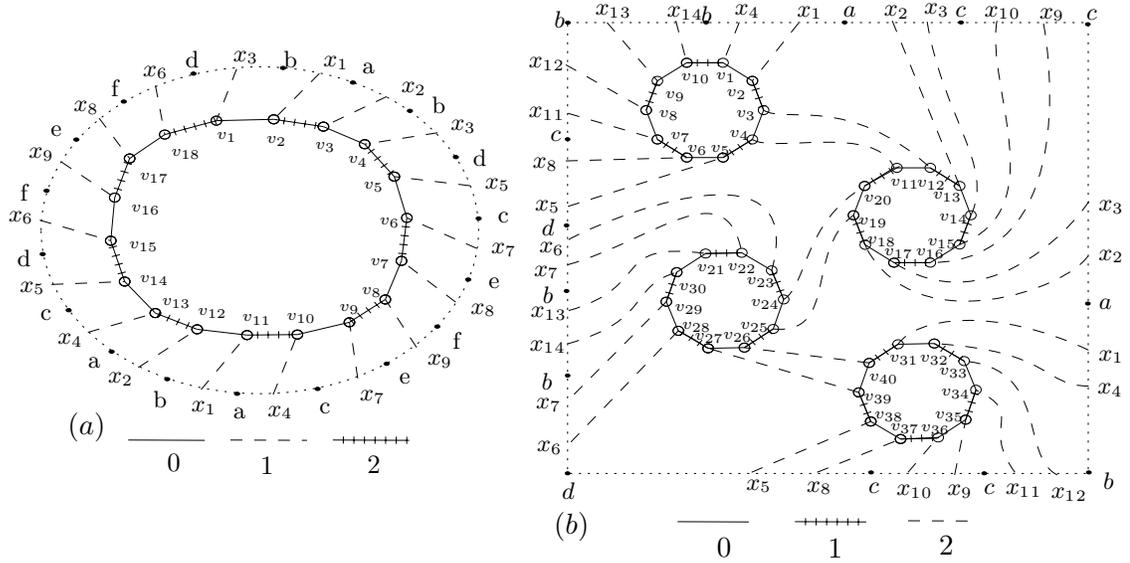

   \begin{center}
\tikzset{every picture/.style={line width=0.3pt}} 

%
\caption{Embedding on $\#_2 \mathbb{T}^2$ of a gem representing $\#_{2} \mathbb{T}^2$: $(a)$  of type $(6^2,18)$; $(b)$ of type $(10^2,4)$.}
    \label{fig:10}

   \end{center}
\end{figure}

\noindent {\bf \underline{$(6^2,12)$-type}}: For the type $(6^2,12)$, the graph has 24 vertices.
 In Figure~\ref{fig:9}$(b)$, the surface consists of six vertices \(a,b,c,d,e,f\), nine edges, and a single face; since the surface is orientable, it is a double torus. Each of the vertices $a$, $d$, $e$, and $f$ lies in the interior of the region bounded by a $(1,2)$)-colored $6$-cycle, while each of the vertices $b$ and $c$ lies in the interior of the region bounded by a $(0,2)$-colored $12$-cycle.
 The \((0,1)\)-colored cycles bound the four inner \(6\)-gons. Consequently, the figure represents a semi-equivelar gem of type \([(6^2,12);24]\) regularly embedded on the surface \(\#_{2}\mathbb{T}^{2}\).

\smallskip

\noindent {\bf \underline{$(6^2,18)$-type}}: For the type $(6^2,18)$, there are 18 vertices.
 In Figure~\ref{fig:10}$(a)$, the surface comprises six vertices \(a,b,c,d,e,f\), nine edges, and a single face; since it is orientable, the surface is a double torus. Each of the vertices \(a,d,e\) lies in the interior of the region enclosed by a \((1,2)\)-colored \(6\)-cycle, while each of the vertices \(b,c,f\) lies in the interior of the region bounded by a \((0,1)\)-colored \(6\)-cycle. The \((0,2)\)-colored cycle bounds the inner \(18\)-gon. Consequently, this configuration yields a semi-equivelar gem of type \([(6^2,18);18]\) regularly embedded on the surface \(\#_{2}\mathbb{T}^{2}\).

\smallskip

\noindent {\bf \underline{$(10^2,4)$-type}}: For the type $(10^2,4)$, the corresponding semi-equivelar graph has 40 vertices.
 In Figure~\ref{fig:10}$(b)$, the surface consists of four vertices \(a,b,c,d\), seven edges, and a single face; since the surface is orientable, it is a double torus. Each of the vertices \(a,b,c,d\) lies in the interior of the region bounded by a \((1,2)\)-colored \(10\)-cycle. The \((0,1)\)-colored cycles bound the four inner \(10\)-gons. Consequently, this figure represents a semi-equivelar gem of type \([(10^2,4);40]\) regularly embedded on the surface \(\#_{2}\mathbb{T}^{2}\).

\begin{figure}[h]
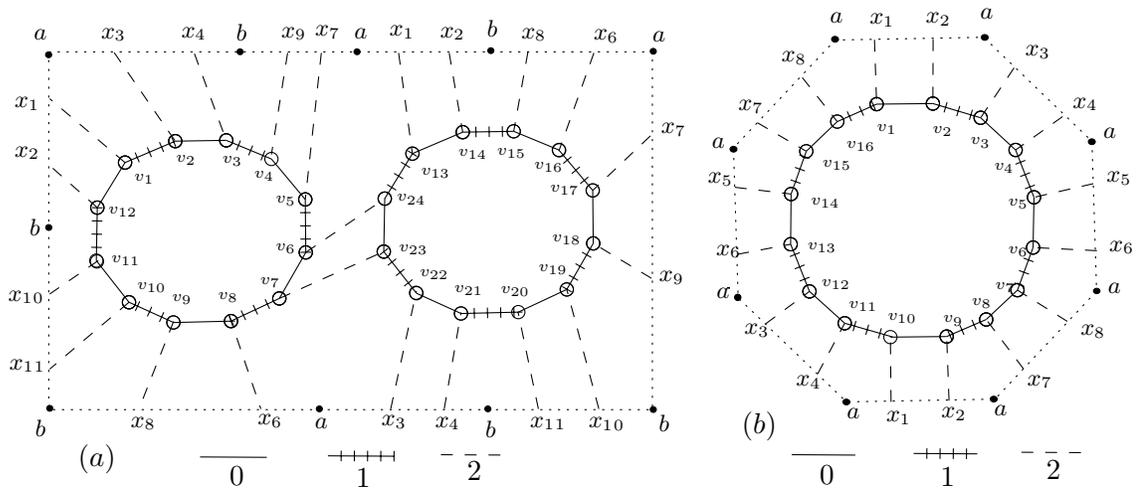

   \begin{center}
\tikzset{every picture/.style={line width=0.4pt}} 


\caption{Embedding on $\#_2 \mathbb{T}^2$ of a gem representing $\#_{2} \mathbb{T}^2$: $(a)$ of type $(12^2,4)$; $(b)$ of type $(16^2,4)$.}
    \label{fig:11}
\end{center}
\end{figure}

\noindent {\bf \underline{$(12^2,4)$-type}}: The type $(12^2,4)$ consists of 24 vertices.
 In Figure~\ref{fig:11}$(a)$, the surface consists of two vertices \(a\) and \(b\), five edges, and a single face; since the surface is orientable, it is a double torus. Each vertex \(a\) and \(b\) lies in the interior of the region bounded by a \((1,2)\)-colored \(12\)-cycle. The \((0,1)\)-colored cycles bound the inner two \(12\)-gons. Consequently, this figure represents a semi-equivelar gem of type \([(12^2,4);24]\) regularly embedded on the surface \(\#_{2}\mathbb{T}^{2}\).

\smallskip

\noindent {\bf \underline{$(16^2,4)$-type}}: The type $(16^2,4)$ consists of 16 vertices.
 In Figure~\ref{fig:11}$(b)$, the surface consists of a single vertex \(a\), four edges, and a single face; since the surface is orientable, it is a double torus. The vertex \(a\) lies in the interior of the region bounded by the \((1,2)\)-colored \(16\)-cycle. The \((0,1)\)-colored cycle bounds the inner \(16\)-gon. Consequently, this figure represents a semi-equivelar gem of type \([(16^2,4);16]\) regularly embedded on the surface \(\#_{2}\mathbb{T}^{2}\).

\smallskip

\noindent {\bf \underline{$(8^2,6)$-type}}: For the type $(8^2,6)$, the corresponding semi-equivelar graph has 24 vertices.
 Figure~\ref{fig:12}$(a)$ depicts a surface with five vertices \(a,b,c,d,e\), eight edges, and a single face; as the surface is orientable, it is a double torus. Each of the vertices \(a,c,e\) is contained in the interior of the region enclosed by a \((0,1)\)-colored \(8\)-cycle, whereas each of the vertices $b$ and $d$ lies inside the region enclosed by a \((1,2)\)-colored \(8\)-cycle. The \((0,2)\)-colored cycles form the boundaries of four inner \(6\)-gons. Hence, the figure corresponds to a semi-equivelar gem of type \([(8^2,6);24]\) regularly embedded on the surface \(\#_{2}\mathbb{T}^{2}\).

 \smallskip

\begin{figure}[h]
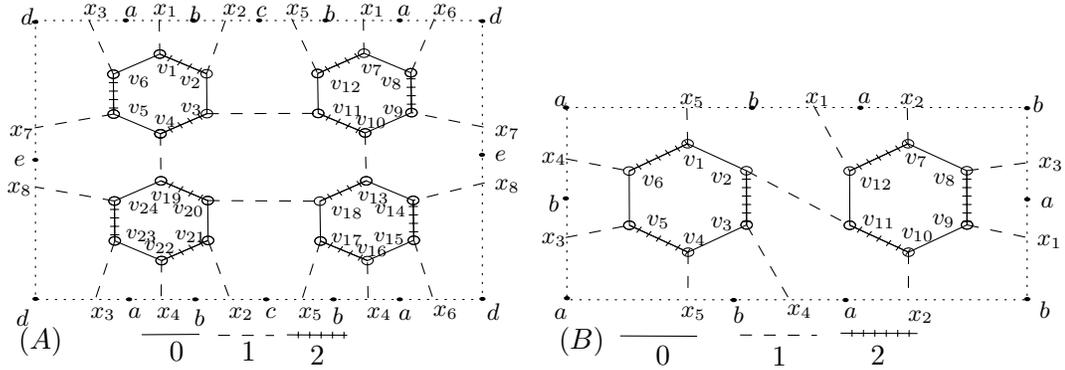

    \begin{center}
\tikzset{every picture/.style={line width=0.4pt}} 



\caption{Embedding on $\#_2 \mathbb{T}^2$ of a gem representing $\#_{2} \mathbb{T}^2$: $(a)$ of type $(8^2,6)$; $(b)$ of type  $(12^2,6)$.}
    \label{fig:12}

 \end{center}    
\end{figure}

\begin{figure}[h]
\begin{center}
\tikzset{every picture/.style={line width=0.7pt}} 

%
\caption{Embedding on $\#_2 \mathbb{T}^2$ of a gem representing $\#_{2} \mathbb{T}^2$ of type $(4,6,14)$.}
    \label{fig:13}
   \end{center}
\end{figure}

\noindent {\bf \underline{$(12^{2},6)$-type}}: For the type $(12^2,6)$, the graph has 12 vertices.
Figure~\ref{fig:12}$(b)$ illustrates a surface with two vertices \(a\) and \(b\), five edges, and a single face; since the surface is orientable, it is homeomorphic to a double torus. The vertex \(a\) is contained in the interior of the region enclosed by the \((1,2)\)-colored \(12\)-cycle, whereas the vertex \(b\) lies in the interior of the region enclosed by the \((0,1)\)-colored \(12\)-cycle. The \((0,2)\)-colored cycles form the boundaries of the two inner \(6\)-gons. Hence, this configuration corresponds to a semi-equivelar gem of type \([(12^{2},6);12]\) regularly embedded on the surface \(\#_{2}\mathbb{T}^{2}\).

\smallskip

\noindent {\bf \underline{$(4,6,14)$-type}}: The type $(4,6,14)$ consists of 168 vertices.
Figure~\ref{fig:13} depicts a surface with sixteen vertices $a,b,c,d,e,f,g,h,i,j,k,l,m,n$, $o,p$, nineteen edges, and a single face; as the surface is orientable, it is a double torus. Each vertex lies in the interior of the region bounded by a $(1,2)$-colored $6$-cycle.
 The \((0,2)\)-colored cycles bound twelve inner \(14\)-gons. Hence, this figure represents a semi-equivelar gem of type \([(4,6,14);168]\) regularly embedded on the surface \(\#_{2}\mathbb{T}^{2}\).

\begin{figure}[h]
\begin{center}
\tikzset{every picture/.style={line width=0.4pt}} 



\caption{Embedding on $\#_2 \mathbb{T}^2$ of a gem representing $\#_{2} \mathbb{T}^2$ of type $(4,6,16)$.}
    \label{fig:14}

\end{center}
\end{figure}

\begin{figure}[h]  
\begin{center}
\tikzset{every picture/.style={line width=0.4pt}} 

%

\caption{Embedding on $\#_2 \mathbb{T}^2$ of a gem representing $\#_{2} \mathbb{T}^2$ of type $(4,6,18)$.}
    \label{fig:15}

\end{center}
\end{figure}

\noindent {\bf \underline{$(4,6,16)$-type}}: The semi-equivelar graph of type $(4,6,16)$ has 96 vertices.
Figure~\ref{fig:14} depicts a surface with four vertices \(a,b,c,d\), seven edges, and a single face; as the surface is orientable, it is a double torus. Each vertex lies inside the region bounded by a \((0,2)\)-colored \(16\)-cycle, while the \((1,2)\)-colored cycles bound sixteen inner \(6\)-gons. Hence, the figure represents a semi-equivelar gem of type \([(4,6,16);96]\) regularly embedded on \(\#_{2}\mathbb{T}^{2}\).

\smallskip

\noindent {\bf \underline{$(4,6,18)$-type}}: The semi-equivelar graph of type $(4,6,18)$ has 72 vertices.
Figure~\ref{fig:15} depicts a surface with four vertices \(a,b,c,d\), seven edges, and a single face; as the surface is orientable, it is a double torus. Each vertex lies inside the region bounded by a \((0,2)\)-colored \(18\)-cycle, while the \((1,2)\)-colored cycles bound twelve inner \(6\)-gons. Hence, the figure represents a semi-equivelar gem of type \([(4,6,18);72]\) regularly embedded on \(\#_{2}\mathbb{T}^{2}\).

\smallskip
\begin{figure}[h]
\begin{center}
\tikzset{every picture/.style={line width=0.4pt}} 


\caption{Embedding on $\#_2 \mathbb{T}^2$ of a gem representing $\#_{2} \mathbb{T}^2$ of type $(4,6,20)$.}
    \label{fig:16}
\end{center}
\end{figure}

\begin{figure}[h]
    \begin{center}
        \tikzset{every picture/.style={line width=0.4pt}} 

%
\caption{Embedding on $\#_2 \mathbb{T}^2$ of a gem representing $\#_{2} \mathbb{T}^2$ of type $(4,6,24)$.}
    \label{fig:17}
    \end{center}
\end{figure}

\noindent {\bf \underline{$(4,6,20)$-type}}: For the type $(4,6,20)$, the graph has 60 vertices.
Figure~\ref{fig:16} illustrates a surface with three vertices \(a,b,c\), six edges, and a single face; since the surface is orientable, it is a double torus. Each vertex lies within the region bounded by a \((0,2)\)-colored \(20\)-cycle, while the \((1,2)\)-colored cycles bound ten inner \(6\)-gons. Thus, the figure represents a semi-equivelar gem of type \([(4,6,20);60]\) regularly embedded on \(\#_{2}\mathbb{T}^{2}\).

\smallskip

\noindent {\bf \underline{$(4,6,24)$-type}}: The type $(4,6,24)$ consists of 48 vertices.
Figure~\ref{fig:17} depicts a surface with two vertices \(a\) and $b$, five edges, and a single face; since the surface is orientable, it is a double torus. Each of the vertices $a$ and $b$ lies inside the region bounded by a \((0,2)\)-colored two \(24\)-cycles, while the \((1,2)\)-colored cycles bound eight inner \(6\)-gons. Hence, the figure represents a semi-equivelar gem of type \([(4,6,24);48]\) regularly embedded on \(\#_{2}\mathbb{T}^{2}\).

\noindent {\bf \underline{$(4,6,36)$-type}}: For the type $(4,6,36)$, the corresponding semi-equivelar graph has 36 vertices.
 In Figure~\ref{fig:18}, the surface consists of a single vertex \(a\), four edges, and a single face; as it is orientable, it is a double torus. The vertex lies within the region bounded by the \((0,2)\)-colored \(36\)-cycle, while the \((1,2)\)-colored cycles bound six inner \(6\)-gons. Thus, the figure represents a semi-equivelar gem of type \([(4,6,36);36]\) regularly embedded on \(\#_{2}\mathbb{T}^{2}\).

\begin{figure}[h]
    \begin{center}

\tikzset{every picture/.style={line width=0.5pt}} 


\caption{Embedding on $\#_2 \mathbb{T}^2$ of a gem representing $\#_{2} \mathbb{T}^2$ of type $(4,6,36)$.}
    \label{fig:18}
    \end{center}
\end{figure}

\noindent {\bf \underline{$(4,8,10)$-type}}: For the type $(4,8,10)$, the graph has 80 vertices.
Figure~\ref{fig:19} depicts a surface with nine vertices \(a,b,c,d,e,f,g,h\), eleven edges, and a single face; since the surface is orientable, it is a double torus. Each of the vertices $a$, $c$, $d$, $e$, $f$, $g$, and $h$ lies in the interior of the region bounded by a $(0,2)$-colored $10$-cycle, while the \((1,2)\)-colored cycles bound ten inner \(8\)-gons. Hence, the figure represents a semi-equivelar gem of type \([(4,8,10);80]\) regularly embedded on \(\#_{2}\mathbb{T}^{2}\).

\begin{figure}[h]
    \begin{center}

\tikzset{every picture/.style={line width=0.4pt}} 



\caption{Embedding on $\#_2 \mathbb{T}^2$ of a gem representing $\#_{2} \mathbb{T}^2$ of type $(4,8,10)$.}
    \label{fig:19}
    \end{center}
\end{figure}

\begin{figure}[h]
    \begin{center}
\tikzset{every picture/.style={line width=0.41pt}} 

%

\caption{Embedding on $\#_2 \mathbb{T}^2$ of a gem representing $\#_{2} \mathbb{T}^2$ of type $(4,8,12)$.}
    \label{fig:20}
    \end{center}
\end{figure}

\noindent {\bf \underline{$(4,8,12)$-type}}: The type $(4,8,12)$ consists of 48 vertices.
Figure~\ref{fig:20} shows a surface with four vertices \(a,b,c,d\), seven edges, and a single face; as it is orientable, it is a double torus. Each vertex lies inside the region bounded by a \((0,2)\)-colored \(12\)-cycle, while the \((1,2)\)-colored cycles enclose six inner \(8\)-gons. Thus, the figure represents a semi-equivelar gem of type \([(4,8,12);48]\) regularly embedded on \(\#_{2}\mathbb{T}^{2}\).

\smallskip

\noindent {\bf \underline{$(4,8,16)$-type}}: The semi-equivelar graph of type $(4,8,16)$ has 32 vertices.
Figure~\ref{fig:21} depicts a surface with two vertices \(a\) and $b$, five edges, and a single face; being orientable, it is a double torus. Each of the vertices $a$ and $b$ lies inside the region bounded by a \((0,2)\)-colored \(16\)-cycle, while the \((1,2)\)-colored cycles enclose four inner \(8\)-gons. Hence, the figure represents a semi-equivelar gem of type \([(4,8,16);32]\) regularly embedded on \(\#_{2}\mathbb{T}^{2}\).

\begin{figure}[h]
    \begin{center}
        
\tikzset{every picture/.style={line width=0.4pt}} 


\caption{Embedding on $\#_2 \mathbb{T}^2$ of a gem representing $\#_{2} \mathbb{T}^2$ of type $(4,8,16)$.}
    \label{fig:21}
    \end{center}
\end{figure}

\noindent {\bf \underline{$(4,8,24)$-type}}: The type $(4,8,24)$ consists of 24 vertices.
Figure~\ref{fig:22}$(a)$ shows a surface with a single vertex \(a\), four edges, and a single face; being orientable, it is a double torus. The vertex lies inside the region bounded by the \((0,2)\)-colored \(24\)-cycle, while the \((1,2)\)-colored cycles enclose three inner \(8\)-gons. Thus, the figure represents a semi-equivelar gem of type \([(4,8,24);24]\) regularly embedded on \(\#_{2}\mathbb{T}^{2}\).

\begin{figure}[h]
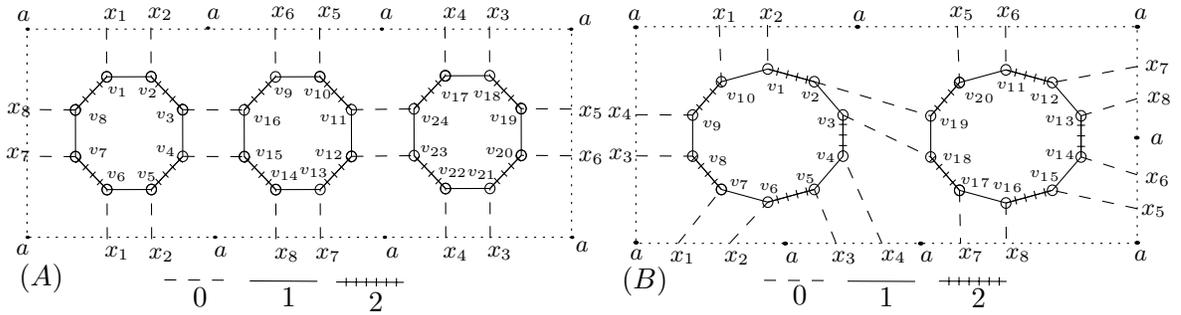

   \begin{center}
       
\tikzset{every picture/.style={line width=0.4pt}} 



\caption{Embedding on $\#_2 \mathbb{T}^2$ of a gem representing $\#_{2} \mathbb{T}^2$: $(a)$ of type $(4,8,24)$; $(b)$ of type $(4,10,20)$.}
    \label{fig:22}
   \end{center}
    \end{figure}

\noindent {\bf \underline{$(4,10,20)$-type}}: The semi-equivelar graph of type $(4,10,20)$ has 20 vertices.
Figure~\ref{fig:22}$(b)$ depicts a surface with a single vertex \(a\), four edges, and a single face; as it is orientable, it is a double torus. The vertex lies within the region bounded by the \((0,2)\)-colored \(20\)-cycle, while the \((1,2)\)-colored cycles enclose two inner \(10\)-gons. Hence, the figure represents a semi-equivelar gem of type \([(4,10,20);20]\) regularly embedded on \(\#_{2}\mathbb{T}^{2}\).
\end{proof}

\section{Conclusion and Future Exploration}\label{conclusion}
In~\cite{ab25, bb24}, the authors classified all semi-equivelar gems of PL $d$-manifolds whose underlying embeddings lie on surfaces of genus $0$ and $1$.  
In the present article, we extend this line of investigation to the case of the double torus, i.e., the orientable surface of genus~$2$.  
The classification problem becomes substantially more complicated for surfaces of genus at least~$3$, as the number of possible configurations grows rapidly, making an exhaustive analysis prohibitively time-consuming.  
Thus, beyond genus~$2$, it is natural to explore alternative directions.

In~\cite{bb24}, we established the uniqueness of semi-equivelar gems of certain types.  
Building on this, one may now investigate characterization problems for each individual type; an interesting direction that may also attract researchers in theoretical computer science due to the combinatorial and algorithmic aspects involved. We therefore propose the following problems.

\medskip

For each $(p_0,p_1,p_2)$-type (where $p_0,p_1,p_2$ need not be distinct), if a gem embeds regularly on $\#_2 \mathbb{T}^2$, then the gem represents the surface $\#_2 \mathbb{T}^2$ itself.  
In this article, we construct examples of regularly embedded gems on $\#_2 \mathbb{T}^2$ for each of the possible $(p_0,p_1,p_2)$-types.  
This leads naturally to the following enumerative question.

\begin{question}
For each possible $(p_0,p_1,p_2)$-type admitting a regular embedding on $\#_2 \mathbb{T}^2$, how many non-isomorphic semi-equivelar gems of that type embed regularly on $\#_2 \mathbb{T}^2$?
\end{question}

We also consider semi-equivelar gems of type $(p_0,p_1,p_2,p_3)$ (where the $p_i$ need not be distinct) that admit regular embeddings in $\#_2 \mathbb{T}^2$.  
For each admissible type, we present an explicit example.  
For the types $(6^4)$, $(4^3,12)$, $(4^2,8^2)$, and $(4,8,4,8)$, the corresponding gems we constructed represent the manifold $\mathbb{S}^3$.  
For the types $(4^3,8)$ and $(4^2,6^2)$, the corresponding gems we constructed represent the manifold $\mathbb{RP}^3$.  The semi-equivelar gems we constructed of types $(4^3,6)$ and $(4,6,4,6)$ represent the lens spaces $L(5,2)$ and $L(3,1)$, respectively.  

Beyond the classification problem for each of these types, one may also consider the following topological realization question.

\begin{question}
For each possible $(p_0,p_1,p_2,p_3)$-type admitting a regular embedding on $\#_2 \mathbb{T}^2$, which non-homeomorphic $3$-manifolds admit a semi-equivelar gem of that type regularly embedded on $\#_2 \mathbb{T}^2$?
\end{question}

In \cite{c92}, the author introduces the notion of locally regular colored graphs as a natural higher-dimensional extension of the classical concept of regular maps on surfaces. Their work develops this framework in dimension three and provides a detailed description of locally regular colored graphs on spherical $3$-manifolds.
Since our results parallel aspects of this theory, a natural direction for future research is to extend the notion of locally regular colored graphs to manifolds that admit a semi-equivelar gem regularly embedded on $\#_2 \mathbb{T}^2$.
Because the gems occurring in our setting exhibit a high degree of symmetry, such an extension would be particularly useful for studying degree-$d$ maps between manifolds that admit semi-equivelar gems of these types. This line of investigation may lead to a broader and more unified framework for understanding symmetries and structural properties of PL manifolds via their colored-graph representations.

The results of this article mark a significant step toward a broader understanding of semi-equivelar gems on higher-genus surfaces. By completing the classification for embeddings on the double torus and constructing explicit examples for all admissible types, we provide a foundational dataset and a framework that will guide future investigations in both combinatorial topology and topological graph theory. The interplay between symmetry, embedding properties, and manifold representation revealed in this work not only deepens the structural insight into PL manifolds but also opens promising avenues for studying locally regular colored graphs and degree-$d$ maps in higher dimensions.

\medskip

\noindent {\bf Acknowledgement:} The authors thank the editor and the anonymous referees for many useful comments and suggestions that have greatly improved the article. The second author is supported by the Mathematical Research Impact Centric Support Research Grant (MTR/2022/000036) by ANRF (India).

\medskip			
				
{\footnotesize

\end{document}